\documentclass{amsart}

\usepackage[T1]{fontenc}
\usepackage{ae,aecompl}

\usepackage{amssymb,amscd}
\usepackage{url}
\usepackage{amsthm,amsfonts,amsmath,amsxtra}
\usepackage{graphicx}
\usepackage{tikz} 

\numberwithin{equation}{section}
\numberwithin{table}{section}
\numberwithin{figure}{section}

\newtheoremstyle{bold}
{.5\baselineskip}{.5\baselineskip}{\itshape}{}{\bfseries}{.}{.5em}{}
\newtheoremstyle{shy}
{.5\baselineskip}{.5\baselineskip}{}{}{\bfseries}{.}{.5em}{}

\makeatletter
\def\@captionfont{\small}

\def\section{\@startsection{section}{1}%
  \z@{.9\linespacing\@plus\linespacing}{.5\linespacing}%
  {\large\bfseries\boldmath\centering}}

\def\subsection{\@startsection{subsection}{2}%
  \z@{.7\linespacing\@plus\linespacing}{.5\linespacing}%
  {\normalfont\scshape\centering}}

\def\theindex{\@restonecoltrue\if@twocolumn\@restonecolfalse\fi
  \columnseprule\z@ \columnsep 35\p@
  \@indextitlestyle
  \thispagestyle{empty}%
  \let\item\@idxitem
  \parindent\z@  \parskip\z@\@plus.3\p@\relax
  \raggedright
  \hyphenpenalty\@M
  \footnotesize}

\renewcommand{\@bibtitlestyle}{%
  \@xp\section\@xp*\@xp{\bibname}%
}
\renewcommand{\tocchapter}[3]{%
  \indentlabel{\@ifnotempty{#2}{\ignorespaces#1 #2.\quad}}#3}

\renewcommand{\tocsection}[3]{%
  \indentlabel{\@ifnotempty{#2}{\makebox[3.2em][l]{\ignorespaces#1 #2.}}}#3}

\makeatother

\renewcommand{\bibname}{References}
\renewcommand{\ge}{\geqslant}

\renewcommand{\le}{\leqslant}

\theoremstyle{bold}
\newtheorem{theorem}{Theorem}[section]
\newtheorem{proposition}[theorem]{Proposition}

\theoremstyle{shy}

\newtheorem{remark}[theorem]{Remark}
\newtheorem{example}[theorem]{Example}

\newcommand{\cA}{\mathcal{A}}

\newcommand{\cL}{\mathcal{L}}

\newcommand{\cP}{\mathcal{P}}

\newcommand{\EE}{\mathbb{E}}

\newcommand{\NN}{\mathbb{N}}

\newcommand{\PP}{\mathbb{P}\ts}

\newcommand{\ts}{\hspace{0.5pt}}

\title[Genealogies and Inference for multiple merger coalescents]{%
 Genealogies and inference for populations with highly
 skewed offspring distributions}
\author{Matthias~Birkner and Jochen~Blath}
\date{\today}
\makeindex

\begin{document}

\begin{abstract}
We review recent progress in the understanding of 
the role of multiple- and simultaneous multiple merger coalescents as models for the genealogy
in idealised and real populations with exceptional reproductive behaviour. In particular, we discuss models with `skewed offspring distribution' (or under other non-classical evolutionary forces) which lead in the single locus haploid case to multiple merger coalescents, and in the multi-locus diploid 
case to simultaneous multiple merger coalescents. Further, we discuss inference methods under the infinitely-many sites model which allow both model
selection and estimation of model parameters under these coalescents.
\end{abstract}

\maketitle

\section{Multiple merger coalescents in population genetics}

\subsection{Introduction}

The `standard' model in 
mathematical population genetics is Kingman's coalescent
\cite{MBJB-K82}, which describes on appropriate time scales the
random genealogies of a large class of population models.
\index{Kingman coalescent}
A salient
feature of models in the domain of attraction of Kingman's coalescent
and its ramifications is that, at least in the limit of large
population size, only {\em binary}\/ mergers of ancestral lineages are
visible.  This is owed to the fact that the number of offspring of any
individual must be negligible in comparison with the total population
size.

It is an important and very useful universality feature of Kingman's
coalescent that as the population size $N\to\infty$, the details of
the actual offspring distribution are `washed out' from the limit
model, only its variance $\sigma_N^2 \to \sigma^2 \in (0, \infty)$
remains as a time-rescaling compared to the `standard' Kingman
coalescent.
A crucial assumption here is $\sigma^2 < \infty$.

The question `what if $\sigma^2=\infty$?' is also biologically
relevant: While all real populations are finite, coalescent
theory is about (tractable) limit results as $N\to\infty$, and
$\sigma^2=\infty$ really means that $\sigma^2_N$ is large when $N$ is
large.
As we will see below, there is a variety of
biological mechanisms which predict a deviation from
the Kingman coalescent model.

In this article, we will first describe 
general coalescent
models (where the term `general' means that 
multiple- and even simultaneous multiple mergers of ancestral
lineages will be allowed), and review briefly population models that lead to limiting
genealogies described by certain subclasses of these general
coalescent processes.  We will then investigate how one of the most
popular statistics of real DNA sequence data (under the infinitely
many sites model), namely the site-frequency spectrum, behaves under
these coalescent models, and then derive inference methods that allow
to estimate evolutionary parameters within a certain class of
coalescent models, or to distinguish between different underlying
genealogical models.  While this theory is mostly confined to
single-locus data of haploid populations, we will finally derive the
genealogy in a simple diploid multi-locus model. Interestingly, this
will naturally lead to genealogies driven by coalescents with
simultaneous multiple mergers.  Also, the additional information
contained in multi-locus data will, despite dependence between
different loci that is inherent in multiple-merger coalescent even in
the face of high recombination rates, increase the statistical power
of our methods for inference.

We conclude this text 
with an outlook on recent developments in the
field and the potential relevance of our results.  To sum up, we aim
to take steps towards understanding in how far the conjecture of Eldon \& Wakeley
(\cite{MBJB-EW06}, p.~2622) holds:
\begin{center} 
 {\center \em `It may be that Kingman's coalescent applies only to a small fraction
of species. For many species, the coalescent with multiple mergers might be
a better null model than Kingman's coalescent.'}
\end{center}
Note that this article is related to several others in this volume that also 
touch upon the topic of non-standard genealogies, in particular those
by Fabian Freund, by G\"{o}tz Kersting and Anton Wakolbinger and by Anja Sturm.
We will highlight concrete links in the sequel.

\subsection{Multiple and simultaneous multiple merger coalescents}

About two decades ago, two natural classes of general coalescent
processes, the so-called $\Lambda$-coalescents
\cite{MBJB-P99,MBJB-S99, MBJB-DK99} and $\Xi$-coalescents
\cite{MBJB-S2000,MBJB-MS2001} were introduced in the mathematical
literature. All these coalescents have in common that they are
(exchangeable) partition-valued continuous-time Markov chains, that
is, they take values in the space
$\cP_n$, the space of finite
partition of $[n] := \{ 1, \ldots, n\}$ if started from a finite
number of blocks.  Both of the above classes of coalescent processes
allow multiple mergers of ancestral lines, by which we mean a
transition that is obtained from the current partition state by
merging a certain number of blocks (representing ancestral lines) into
one or several new blocks, thus obtaining a `coarser partition'. In
the case of the classical Kingman coalescent, these transitions are always
binary, that is, precisely two blocks merge into one new block.

In the case of a $\Lambda$-coalescent, however, at transition times,
multiple lines necessarily merge into one single new block, while for
$\Xi$-coalescents, subsets of blocks involved in a coalescence event
may merge into different `target blocks'.

The path of an $n$-coalescent process corresponds in a natural way to a random tree
where the leaves correspond to $\{1\}, \{2\}, \dots, \{n\}$ and internal nodes
to larger blocks. 
In fact, one can interpret a coalescent as a random
metric space; see e.g.\ \cite{MBJB-GPW09} and \cite{MBJB-Gufler2018a, MBJB-Gufler2018b}.

In this article, we only consider coalescent processes starting from finitely many blocks 
(i.e., $n$-coalescents). The corresponding coalescents with $n=\infty$ can be constructed by employing 
consistency and using Kolmogorov's extension theorem, or explicitly via look-down constructions 
\cite{MBJB-DK99, MBJB-BBM09}. They have very interesting mathematical properties which are, however, 
not in the focus of this text. Let us first briefly introduce the pertinent notation. 
\smallskip

\subsubsection{Multiple merger (MMC) coalescents}

\index{$\Lambda$-coalescent}
\index{multiple merger coalescent}
\index{coalescent!multiple merger}
For $\pi \in \cP_n$ let $|\pi|$ denote the number of blocks and for 
$\pi, \pi' \in \cP_n$ we write $\pi^\prime \prec_{m,k} \pi$ if $|\pi|=m$ and $\pi'$
arises from $\pi$ by merging $k$ blocks into a single one (a `$k$-merger'). 

For a finite measure $\Lambda$ on $[0,1]$, define 
\begin{equation}
\label{MBJB-eq:lambda_m_k}
\lambda_{m,k}:=\int_0^1 x^{k-2}(1-x)^{m-k}\Lambda(dx), \quad \lambda_{m} := \sum_{k=2}^m
\binom{m}{k} \int_0^1 x^{k-2}(1-x)^{m-k}\Lambda(dx).
\end{equation} 
The $n$-$\Lambda$-coalescent is a $\cP_n$-valued continuous-time Markov chain  $\{\Pi_t^{(\Lambda)},t\ge0\}$ 
with  transition rates   $q_{\pi, \pi^\prime}$  from  $\pi$ to $\pi^\prime \neq \pi$ given by 
\begin{equation}
\label{MBJB-eq:lambdarates}
    q_{\pi, \pi^\prime} = \begin{cases} \lambda_{m,k}
      & \textrm{if $ \pi^\prime \prec_{m,k} \pi$ for some $k$,}\\
    0 & \textrm{otherwise.} \\ \end{cases}
\end{equation}
\begin{remark}
  \label{MBJB-rem:lambda_m_k interpretation}
  A natural interpretation of \eqref{MBJB-eq:lambda_m_k} is to imagine that for $x \in (0,1]$ at rate
  $x^{-2} \Lambda(dx)$, a `merging event of size $x$' occurs: In such an event, every
  block independently flips a `coin' with success probability $x$ and all the `successful'
  blocks are merged. In fact, such constructions are in \cite{MBJB-P99, MBJB-DK99}
  and this intuition is also corroborated by the duality with the $\Lambda$-Fleming-Viot
  process 
  (see page~\pageref{MBJB-pageref:Lambda-FV}).
\end{remark}
\smallskip

Obviously, the class of all $\Lambda$-coalescents (corresponding to
all the finite measures on $[0,1]$) is quite large and in particular
non-parametric.  The following important special cases have frequently
appeared in the literature:
\begin{example}
  \label{MBJB-ex:Lambdacoals}
  \begin{itemize}
  \item[{\tt (K)}] The Kingman coalescent $\Pi^{\tt (K)}$ \cite{MBJB-K82}
     corresponds to the choice
    \begin{equation*}
      \label{eq:Kingman}
      \Lambda(dx)= \delta_0(dx).
    \end{equation*} 
    i.e.\ $\Pi^{\tt (K)} = \Pi^{(\delta_0)}$. 
    Here, the measure
    $\Lambda$ is concentrated on the point $0$ and no multiple, only
    binary mergers happen, as is evident from
    \eqref{MBJB-eq:lambda_m_k}.
    \index{coalescent!Kingman}
    \index{Kingman coalescent}
  \item[{\tt (S)}] The `star-shaped coalescent' coalescent $\Pi^{\tt (S)}$ corresponds to the choice
    \begin{equation*}
      \label{eq:star-shaed}
      \Lambda(dx)= \delta_1(dx).
    \end{equation*} 
    This coalescent exhibits only one single transition, in which all active lines merge into a single line within one step.
  \item[{\tt (BS)}] The Bolthausen-Sznitman coalescent $\Pi^{\tt (BS)}$, introduced in \cite{MBJB-BS98}
    as a tool to study certain spin glass models in statistical mechanics, is given by 
    \begin{equation*}
      \label{eq:b-s}
      \Lambda(dx)= {\bf 1}_{[0,1]}(x)(dx),
    \end{equation*} 
    i.e.\ when the measure $\Lambda$ is the uniform distribution on $[0,1]$.  
  \item[{\tt (B)}] The Beta$(2-\alpha, \alpha)$-coalescent $\Pi^{\tt (B)}$ is given by
    \begin{equation*}
      \label{eq:betameasure} 
      \Lambda(dx) = \frac{\Gamma(2)}{\Gamma(2-\alpha)\Gamma(\alpha)} x^{1-\alpha} (1-x)^{\alpha-1} \,dx,
    \end{equation*}
    with $\alpha \in (0,2)$. Here, the measure $\Lambda$ is associated
    with the beta distribution with parameters $2-\alpha$ and $\alpha$.
    The limiting case $\alpha=2$ (in the sense of weak convergence of
    measures) corresponds to the Kingman coalescent, while $\alpha=1$ 
    returns the Bolthausen-Sznitman-coalescent $\Pi^{\tt (BS)}$ and (the weak limit) $\alpha \to 0$ gives the star-shaped coalescent $\Pi^{\tt (S)}$.

    For a visual impression of realisations of Beta-coalescent trees for different values of $\alpha$ 
    we refer to the contribution by G\"{o}tz Kersting and Anton Wakolbinger in this volume.
    in the article by G.~Kersting and A.~Wakolbinger in
    this volume.
    \index{Beta-coalescent}
  \item[{\tt (EW)}] The following class of purely atomic coalescents has been investigated by \cite{MBJB-EW06}: Here, one considers the cases
    \begin{equation*}
      \label{eq:EWoneatom}
      \Lambda(dx) =  \delta_\psi(dx), 
    \end{equation*}
    and
    \begin{equation*}
      \label{eq:EWtwoatom}
      \Lambda(dx) =\frac{2}{2+\psi^2} \delta_0(dx) + \frac{\psi^2}{2+\psi^2} \delta_\psi(dx), 
    \end{equation*}  
    with $\psi \in [0,1]$, where $\psi=0$ gives the Kingman coalescent.
    \index{Eldon-Wakeley coalescent}
  \end{itemize}
\end{example}

We refer to \cite{MBJB-GIM14} and \cite{MBJB-Ber09} for surveys on $\Lambda$-coalescents.
See also the contribution by G.~Kersting and A.~Wakolbinger in this volume.
\smallskip

\subsubsection{Simultaneous multiple merger (SMMC) coalescents}
\index{$\Xi$-coalescent}
\index{Simultaneous multiple merger (SMMC) coalescent}
\index{coalescent!simultaneous multiple merger}

Formulating the dynamics of a SMMC requires some notational overhead but we will see that they appear
naturally as genealogies in diploid population models with highly skewed offspring distributions.
For 
\begin{equation}
  \label{MBJB-eq:kgroupmerger}
  \underline{k}=(k_1, k_2, \dots, k_r) \quad \text{with } r \in \NN, \; k_1 \ge k_2 \ge \cdots \ge k_r \ge 2
\end{equation}
and $\pi, \pi' \in \cP_n$ with $|\pi|=m$ we write $\pi' \prec_{m,\underline{k}} \pi$ if $\pi'$
arises from $\pi$ by merging $r$ groups of blocks of sizes $k_1, k_2, \dots, k_r$ 
(and leaving the other blocks unchanged). We write $|\underline{k}|=k_1+\cdots+k_r$.

In order to describe the dynamics of a SMMC, we need a bit of notation: 
Let $\Delta$ denote the infinite simplex 
\begin{equation*}
  \label{eq:simplex} 
  \Delta := \left\{ \ensuremath{\boldsymbol{x}} = (x_1,
    x_2, \ldots ) : x_1\ge x_2\ge\cdots\ge0,\quad\sum_i x_i \le 1
  \right\}
\end{equation*}
and let
$\Delta_{\ensuremath{\boldsymbol{0}}} :=  \Delta\setminus \{(0,0, \ldots ) \} = \Delta\setminus \{\ensuremath{\boldsymbol{0}}\}.$
Let $\Xi_0$ be a finite measure on $\Delta_{\ensuremath{\boldsymbol{0}}}$, $a>0$, then $\Xi:= a \delta_{\mathbf{0}} + \Xi_0$ 
is a finite measure on $\Delta$. 

For $\underline{k}$ as in \eqref{MBJB-eq:kgroupmerger}, with 
$s=m-|\underline{k}|$, put 
\begin{align}
  \label{MBJB-eq:Xilambda_mk}
  \lambda_{m,\underline{k}} & = 
                              a\ensuremath{\boldsymbol{1}}_{(r=1, k_1 = 2)} \notag \\
                            & \quad + \int\limits_{\Delta_{\ensuremath{\boldsymbol{0}}}} \dfrac{ \sum\limits_{\ell = 0}^s \sum\limits_{i_1 \neq \ldots \neq i_{r+\ell }}\binom{s}{\ell}x_{i_1}^{k_1}\cdots x_{i_r}^{k_r} x_{i_{r+1}}\cdots  x_{i_{r   + \ell }} \left( 1 - \sum_j x_j  \right)^{s-\ell}   }{\sum_{j} x_j^2 } \,\Xi_{\ensuremath{\boldsymbol{0}}} (d\ensuremath{\boldsymbol{x}})
\end{align}

An $n$-$\Xi$-coalescent $\{ \Pi^\Xi_t \}$ is a continuous-time 
Markov chain on $\cP_n$ which jumps from $\pi \in \cP_n$ with $|\pi|=m$ to 
$\pi' \in \cP_n$ at rate $q_{\pi, \pi'} = \lambda_{m,\underline{k}}$ if $\pi' \prec_{m,\underline{k}}$ with 
$\underline{k}$ as in \eqref{MBJB-eq:kgroupmerger}, and $q_{\pi, \pi'}=0$ if $\pi' \neq \pi$ is not of this 
form.

The form of the jump rates \eqref{MBJB-eq:Xilambda_mk} has a similar interpretation as 
discussed in Remark~\ref{MBJB-rem:lambda_m_k interpretation} for the case of $\Lambda$-coalescents: 
At rate $a$, pairwise merging occurs. Furthermore, for $\mathbf{x} = (x_1,x_2,\dots) \in \Delta_{\ensuremath{\boldsymbol{0}}}$, 
at rate $(\sum_j x_j^2)^{-1} \Xi_0(dx)$ an `$x$-merging event' occurs. In such an event, every block independently 
draws a `colour,' where colour $i$ is drawn with probability $x_i$ for $i \ge 1$ and colour $0$ with 
probability $1-|\mathbf{x}|$. Then all blocks with the same colour $i$ for $i \ge 1$ are merged. 

\smallskip

Obviously, the class of $\Xi$-coalescents is even richer than the
class of $\Lambda$-coalescents. In particular, one recovers a
$\Lambda$-coalescent by choosing $\Xi :=  \Lambda\otimes \delta_0 \otimes \delta_0 \otimes \cdots,$
i.e.\ if $\Xi$ is concentrated on the first component of the simplex.
However, only a handful of natural
examples have been motivated and analysed on the basis of an
underlying population model so far. The following important special
cases have appeared in the literature:
\begin{example}
  \label{MBJB-ex:Xicoals}
  \begin{itemize}
  \item[{\tt (PD)}]
    Let $PD_\theta$ be the Poisson-Dirichlet distribution with $\theta>0$.
    The Poisson-Dirichlet coalescent with $\Xi = \left( \sum_i x_i^2 \right)^{-1} \! PD_\theta$
    appears in \cite{MBJB-S2003} as the genealogy of the `Dirichlet compound Wright--Fisher model.'
    \index{Poisson-Dirichlet coalescent}
  \item[{\tt (SK)}] Subordinated Kingman-coalescents. If one applies a discontinuous time-change to a Kingman coalescent, as soon as more than one binary coalescence event of the original process falls into a jump-interval of the time-change, one obtains a multiple or simultaneous multiple merger event.
    When the (random) time-change is given by a subordinator $\{S_t\}$, the time-changed process
    $
    \{\Pi^{\tt (K)}_{S_t}\}_{t \ge 0}
    $
    is a $\Xi$-coalescent. The representation of $\Xi$ in terms of $\{S_t\}$ as mixture of Dirichlet distributions is non-trivial and omitted here for brevity, see \cite[Prop.~6.3]{MBJB-BBM09} for a partial answer. See also \cite{MBJB-GMS19} for the related class of `symmetric coalescents'.
       \index{Kingman coalescent!subordinated}
  \item[{\tt (DS)}] R.~Durrett and J.~Schweinsberg \cite{MBJB-DS05} approximate the
    genealogy in a selective sweep by a $\Xi$-coalescent, where $\Xi$ is described by
    a stick-breaking construction, see \cite[Section~3]{MBJB-DS05}.

  \item[{\tt (xEW), (xB)}] In diploid bi-parental populations, in which the reproduction events of each parent are governed by a certain $\Lambda$-coalescent, one obtains genealogies given by $\Xi-$coalescents of the form 
    $$
    \Xi= \frac 14 \int_{[0,1]} \delta_{(x/4,\, x/4,\, x/4,\, x/4,\, 0,\, 0,\, 0,\, \dots)} \, \Lambda(dx) 
    $$
    In particular, the cases $\Lambda = \delta_\psi$ and $\Lambda = \mathrm{Beta}(2-\alpha, \alpha)$ for suitable $\psi$ and $\alpha$ have been considered, see \cite{MBJB-BBE13}. The reason for the fourfold split is that the ancestral line of a chromosome may merge into any of the four parental chromosome (two for each parent).
    Such $\Xi$-coalescents will play an important role in Section~\ref{MBJB-sect:multilocus etc} below.
    \index{diploid}
  \end{itemize}
\end{example}

\subsection{Population models}
\label{MBJB-sect:popmodels}
\index{population model}

A substantial amount of work has been devoted to understanding conditions under which population models converge to limits whose genealogy can be described by one of the above coalescent processes. Typically, one considers populations of fixed size $N$, whose reproductive event 
can be described by {\em exchangeable} offspring distributions. 

A full classification of offspring distributions and time scalings in Cannings-models for convergence to $\Lambda$- and $\Xi$-coalescents has been found in \cite{MBJB-MS2001}. It is thus possible to provide abstract criteria and descriptions for population models that make their ancestral distributions converge to 
any prespecified $\Xi-$ or $\Lambda-$coalescent. 
\index{Cannings model}

However, the relevance of a particular (SMMC) model clearly depends on its plausibility as limit of a in some sense {\em natural} population model. We thus now briefly review such population models and their genealogical coalescent limits.

\begin{itemize}
\item[{\tt (B)}] Beta$(2-\alpha,\alpha)$-coalescents with
  $\alpha \in (1,2]$ are obtained as limiting genealogy of
  Schweinsberg's model \cite{MBJB-S03}, in which individuals produce
  in a first step potential offspring according to a stable law with
  index $\alpha$ and mean $m>1$, and then $N$ out of these are selected for
  survival. This corresponds to what is known as a `highly skewed
  offspring distribution' or `sweepstakes reproduction' (cf.\
  \cite{MBJB-A04, MBJB-H94, MBJB-HP11}). In population biology, it
  resembles so-called `type-III survivorship', that is, high fertility
  leading to excessive amounts of offspring, corresponding to the
  first reproduction step, whereas high mortality early in life is
  modelled in the second step. Several authors have proposed this class
  of coalescents to describe the reproductive behaviour of 
    Atlantic cod (see e.g.\ \cite{MBJB-SBB13, MBJB-AH2015}).   

  One can see heuristically why this particular form of the $\Lambda$-measure 
  appears: The probability that a given individual's offspring provides more than 
  fraction $y$ of the next generation, given that the family is 
  substantial (i.e.\ given $X_1 \ge \varepsilon N$, for $y > \varepsilon$), is 
  approximately
  \begingroup
    \allowdisplaybreaks
  \begin{align*}
    \hspace{4em} \PP\Big( & \frac{X_1}{X_1 + (N-1) m} \ge y
    \mid X_1 \ge \varepsilon N \Big)\\
    & \hspace{4em} = \PP\Big( {X_1} \ge \frac{(N-1) m y}{1-y} \mid X_1 \ge \varepsilon N \Big) \\
    & \hspace{4em} \sim \, \text{const.} \times \frac{(1-y)^\alpha}{y^\alpha}
      = \text{const.} \times \mathrm{Beta}(2-\alpha, \alpha)([y,1]),
  \end{align*}
  \endgroup
    where we replaced $X_2+\cdots+X_n \approx (N-1)m$ by the law of large numbers.
  The model is also mathematically appealing, since it exhibits a close
  connection to renormalised $\alpha$-stable branching processes, see
  \cite{MBJB-BBC05+}.
  \index{Schweinsberg's model}
\item[\texttt{(B')}]
  Huillet's Pareto model: \cite{MBJB-H14} derives $\mathrm{Beta}(2-\alpha,\alpha)$-coalescents
  as limiting genealogies in a population model similar to the one in \texttt{(B)}
  where the sampling can be interpreted as according to a `random fitness value.'
 
\item[{\tt (BS)}] The Bolthausen-Sznitman coalescent appears for
  $\alpha=1$ in the sweepstakes model, but also as limiting genealogy
  at the `tip of a fitness wave.' This was predicted in \cite{MBJB-DWF13} using non-rigorous arguments
  (for a related model also \cite{MBJB-NH13}),
  and partly confirmed (for certain variations of the model) in \cite{MBJB-BBS13}, \cite{MBJB-S17a, MBJB-S17b}.

\item[{\tt (EW)}] This model corresponds to populations, in which in each reproductive step, a fraction of $\psi$ individuals are produced by one single parent. This can be combined with classical Wright-Fisher type reproduction to produce the `Kingman atom' at 0. See \cite{MBJB-EW06}.

\item[{\tt (GM)}] Generalised Moran models. Independently in each reproduction event, a random number $\Psi^{(N)}$
  of offspring are born to a single pair of parents, these offspring replace $\Psi^{(N)}$ randomly chosen
  individuals from the present population. $\PP(\Psi^{(N)}=1)=1$ corresponds to the classical
  Moran model; \texttt{(EW)} is also a special case of this. 
  By suitably choosing $\cL(\Psi^{(N)})$ one can in fact approximate any $\Lambda$-coalescent, see 
  Section~\ref{MBJB-subsect:diploidmodel}. 

\item[{\tt (xEW), (xB)}] Appear as scaling limits of diploid bi-parental models with skewed reproduction. We will present a corresponding model in Section~\ref{MBJB-subsect:diploidmodel}. 
  A complete classification of the corresponding diploid population limits can be found in \cite{MBJB-BLS18}.
 
\end{itemize}

See also Tellier and Lemaire \cite{MBJB-TL14} for a recent overview from a biological perspective.
There are many further extensions of population and coalescent models in the literature, including spatial models such as Barton, Etheridge and V\'eber's spatial $\Lambda$-Fleming Viot process \cite{MBJB-BEV10}, or so-called on/off coalescents in situations with seed banks, see, e.g., the contribution by the second author together with Noemi Kurt in this volume.
However, in this article, our focus is the reproductive mechanism of neutral well-mixed populations, so that we refrain from providing a further discussion of these models here.
\smallskip

All of the above coalescent processes are {\em dual} to the corresponding forward-in-time population limit, given as a (generalised) Fleming-Viot process (which is a measure-valued (jump-)diffusion),  \label{MBJB-pageref:Lambda-FV}
\cite{MBJB-DK99} and e.g.\ \cite{MBJB-BLG03}. 

Details of this and a representation of the generator of $\Xi$-coalescents can be found in \cite{MBJB-BBM09}. There, it is also shown that the above duality can be strengthened to a strong pathwise duality via an extension of Donnelly and Kurtz' celebrated {\em lookdown-construction} \cite{MBJB-DK96, MBJB-DK99}.

\section{Inference based on the site-frequency spectrum}

\index{sequence data}
\index{site-frequency spectrum}
\index{infinitely-many-sites model}
One of the most important and well-studied statistical quantities
derived from DNA sequence data is the {\em site frequency spectrum}
(SFS)\footnote{One can in fact attempt to base statistical inference on the likelihood 
of the full sequence data, see e.g.\ \cite{MBJB-SBB13} and references there. However, this is computationally still prohibitively 
expensive even for moderate sample sizes.}. For the theoretical analysis, we assume that all underlying data fits to the {\em infinitely-many-sites model} (IMS) of population genetics (cf.\ \cite{MBJB-W75} 
or \cite{MBJB-W08}), that is, we assume that every observed site
mutated at most once during the entire history of the sample. This assumption is often at least approximately
true since typical per-site mutation rates are very small.
Here, `site' refers to a single base pair in the DNA molecule. 
Furthermore, from a pragmatical point of view,
the SFS of a dataset is well-defined even if the assumptions of the IMS model are violated 
(see, e.g., \cite{MBJB-Gusfield91} for the combinatorial characterisation of data complying with the IMS model).

\begin{figure}[h!]
  \includegraphics[width=6.75cm, height=3cm]{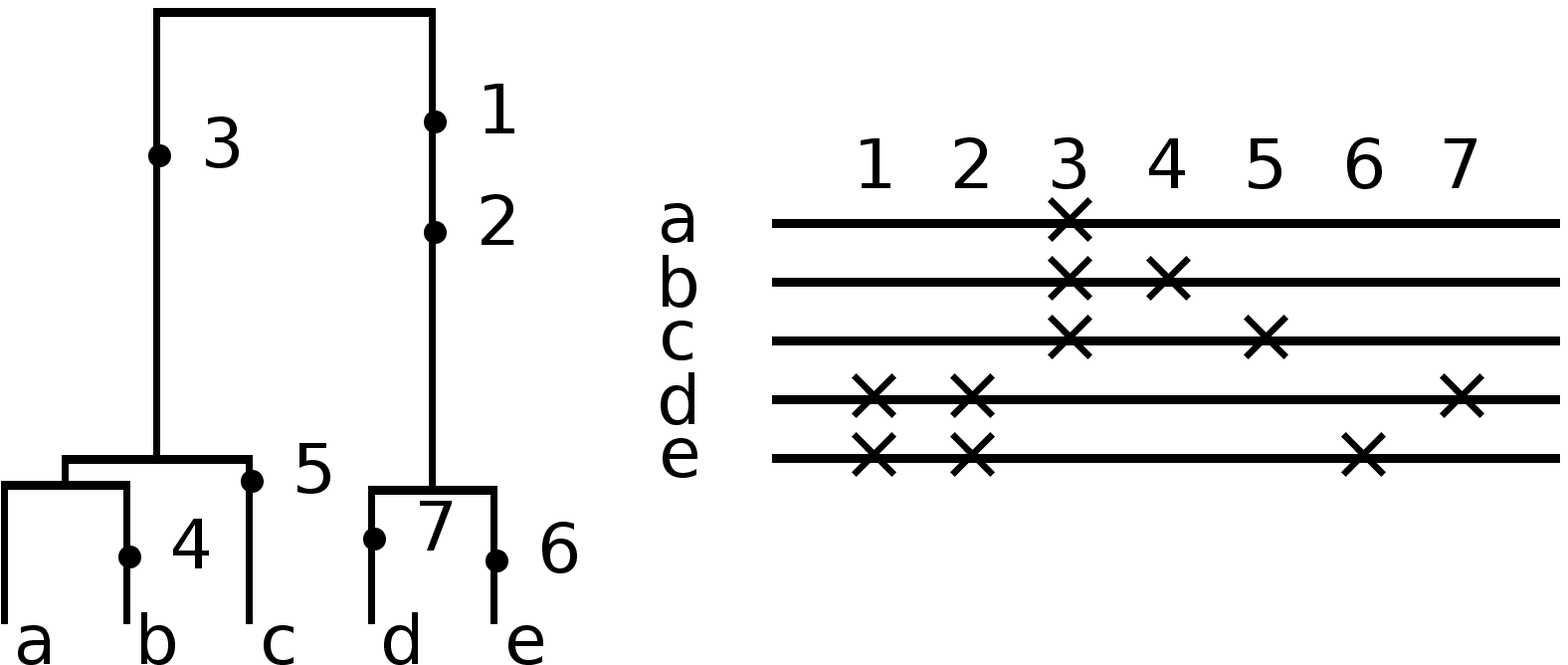}
   \hspace{0.8cm}%
   \parbox[b]{4cm}{\small Mutations on a 
     coalescent tree and resulting data matrix 
     (in schematic 
     form). Implicitly, identical columns are removed from the data matrix.
     The corresponding SFS is ${\boldsymbol{\xi}}^{(5)}=(4,2,1,0)$.}
   \caption{\ }
   \label{MBJB-fig:coalIMS}
 \end{figure}

For the analysis, we also assume that the genealogy of a sample of size $n\in \NN$ is described by one of the above coalescent models $\Pi$ and that mutations occur at some rate $\theta/2>0$ on the coalescent
branches, see Figure~\ref{MBJB-fig:coalIMS} for an illustration. If we know the ancestral state, then, the SFS of an
$n$-sample is defined as
$$
{\boldsymbol{\xi}}^{(n)} := \big(\xi_1^{(n)}, \dots, \xi^{(n)}_{n-1}\big),
$$
where $\xi^{(n)}_i, i \in [n-1]$ is the number of sites at
which a mutation appears $i$-times in our sample.

If the ancestral states are unknown (and thus the data matrix as in Figure~\ref{MBJB-fig:coalIMS} is only defined up to
column-flips), one considers instead the {\em folded} site frequency
spectrum ($\delta_{i, j}$ is the Kronecker delta)
$$
{\boldsymbol{\eta}}^{(n)}:= \Big(\eta_1^{(n)}, \dots, \eta^{(n)}_{\lfloor n/2 \rfloor}\Big)
\quad \text{with } \eta^{(n)}_i = 
\xi^{(n)}_i + \big( 1-\delta_{i, n-i}\big) \xi^{(n)}_{n-i},
\quad i=1, \dots, \lfloor n/2 \rfloor
.$$
\index{site-frequency spectrum!folded}

\subsection{The expected site frequency spectrum}

For a coalescent process $\Pi=\{\Pi_t\}_{t\ge 0}$ with mutation rate $\theta$ we denote its law by $\PP^{\Pi, \theta}$, that is, the law of the coalescent process $\Pi$ on which mutations appear along its branches at rate $\theta/2$. We denote the expectation corresponding to $\PP^{\Pi, \theta}$ by $\EE^{\Pi, \theta}$.
\index{block-counting process of a coalescent}
Recall that the {\em block-counting process} $Y = \{ Y_t \}_{t \ge 0}$ of the  coalescent process $\Pi$ 
\begin{equation}
\label{MBJB-eq:BCP}
Y_t := |\Pi_t|, \quad t \ge 0,
\end{equation}
simply counts the number of ancestral lineages present at each time. 
Then, a general representation of $\EE^{\Pi, \theta}\big[\xi_i^{(n)}\big]$ for any coalescent model $\Pi$
(see \cite{MBJB-GT98}) is 
\begin{equation}
\label{MBJB-eq:masterSFS}
\EE^{\Pi, \theta}\left[\xi_i^{(n)}\right] 
=\frac{\theta}{2}\sum_{k=2}^{n-i+1} {p^{(n), \Pi}[k,i]} \cdot k \cdot \EE^{\Pi}\left[ T_k^{(n)}\right], \quad i \in [n-1],
\end{equation}
where $T_k^{(n)}$ is the random amount of time that $\{ Y_t \}_{t \ge 0}$, starting from
$Y_0=n$, spends in state $k$, and 
$p^{(n), \Pi}[k,i]$ is the probability that \emph{conditional} on the event that $Y_t=k$ for some time point $t$, 
a given one of these $k$ blocks subtends exactly $i \in [n-1]$ leaves.
Thus, in \eqref{MBJB-eq:masterSFS} mutations are classified according to the `level' $k$, which is the value of the block-counting process when they appear in the tree. 

\smallskip

\subsubsection{The block-counting process}
For brevity, we consider only $\Lambda$-coalescents $\Pi$ in this paragraph. 
We see from \eqref{MBJB-eq:lambdarates} 
that $Y$ corresponding to $\Pi$ from \eqref{MBJB-eq:BCP} is itself a continuous-time Markov chain on $\NN$ 
(as $\lambda_{\pi,\pi'}$ depends only on $\pi$ and $\pi'$)
with jump rates
\[ 
q_{ij} = {i \choose i-j+1} \lambda_{i, i-j+1}, \quad i > j \ge 1.
\]
The total jump rate away from state $i$ is $-q_{ii} = \sum_{j=1}^{i-1} q_{ij}$.

We will need the Green function of $Y$, 
\begin{equation}%
\label{MBJB-eq:YtGreenfct}
g(n,m) := \EE_n\left[ \int_0^\infty {\bf 1}_{( Y_s = m )} \, ds \right] 
\quad \mbox{for $\quad n \ge m \ge 2$} .
\end{equation}
For the Kingman coalescent, we have $g(n,m) = \frac{2}{m(m-1)}$ for $m \le n$, 
for the Bolthausen-Sznitman coalescent, explicit expressions can be obtained from 
\cite{MBJB-MP14}. 
In general, there is no explicit formula for \eqref{MBJB-eq:YtGreenfct},  
but decomposing according to the first jump of $Y$ gives a recursion for $g(n,m)$: 
\begin{equation}
\label{MBJB-eq:grec} 
g(n,m)  =  \sum_{k=m}^{n-1} p_{nk} g(k,m), \;\;
n > m \ge 2, \quad \mbox{ and } \quad
g(m,m) =  \frac{1}{-q_{mm}}, \;\; m \ge 2  
\end{equation}
where 
$ p_{nk} := \frac{q_{nk}}{-q_{nn}}$
are the transition probabilities of the embedded discrete skeleton chain.

\subsubsection{\bf The expected SFS for $\Lambda$-coalescents}

\index{site-frequency spectrum!expected}
Decomposing according to the first jump of $Y$ corresponding to a $\Lambda$-coalescent $\Pi$, starting from $n$, 
yields a recursion for $p^{(n),\Lambda}[k,b]$: 
\begin{proposition}[{\cite[Proposition~1 and Proposition~A.1]{MBJB-BBE13}}] 
\label{MBJB-prop:pLambda_kb}
For $1 < k \le n$, we have
\begin{align}
\label{MBJB-eq:pLambda_kb}
{p^{(n), \Lambda}[k,b]} = 
\sum_{n'=k}^{n-1} p_{n,n'}\frac{g(n',k)}{g(n,k)} 
\bigg( &{\bf 1}_{(b > n-n')} \frac{b-(n-n')}{n'} {p^{(n'), \Lambda}[k,b-(n-n')]}\\
&+ {\bf 1}_{(b < n')} \frac{n'-b}{n'} {p^{(n'), \Lambda}[k,b]} \bigg), \notag
\end{align}
with the boundary conditions ${p^{(n), \Lambda}[n,b]}=\delta_{1b}$ and 
${p^{(n), \Lambda}[k,b]} = 0$ if $b > n-(k-1)$.
\end{proposition}
The terms on the right-hand side of \eqref{MBJB-eq:pLambda_kb}
have a natural 
interpretation:
The probability of seeing a jump from $n$ to $n'$, conditionally on
hitting $k$, has probability $p_{n,n'}\frac{g(n',k)}{g(n,k)}$.
Namely, by the Markov property of $Y$,
\[
  \frac{\PP_n \{Y \mbox{ first jumps
      to } n' \cap Y \mbox{ hits } k \}} {\PP_n \{ Y \mbox{ hits } k
    \}} = p_{n, n'} \frac{\PP_{n'} \{ Y \mbox{ hits } k \}} {\PP_n \{ Y
    \mbox{ hits } k \}} = p_{n, n'} \frac{g(n', k)}{g(n, k)}.
\]
Then, thinking `forwards in time from $n'$ lineages', either the
initial $(n-n'+1)$-split occurred to one of the (then necessarily
$b-(n-n')\,$) lineages subtended to the one we are interested in, or
it occurs to one of the (then necessarily $n'-b$) others.  \medskip

Specialising \eqref{MBJB-eq:masterSFS} to the case of a $\Lambda$-coalescent $\Pi$, 
combined with $\EE^{\Pi}\left[ T_k^{(n)}\right] = g(n,k)$ 
(with $g(n,k)$ from \eqref{MBJB-eq:YtGreenfct}, which can be computed recursively via \eqref{MBJB-eq:grec}) gives 
\begin{proposition}
\label{MBJB-prop:Exi}
We have, for $i = 1, \dots, n-1$,
\begin{equation}
  \label{MBJB-eq:Exii}
\EE^{\Lambda, \theta}\left[\xi_i^{(n)}\right] = \frac{\theta}{2}\sum_{k=2}^{n-i+1} {p^{(n), \Lambda}[k,i]} \cdot k \cdot g(n,k).
\end{equation}
\end{proposition}

\begin{figure} 
\includegraphics[width=5.5cm]{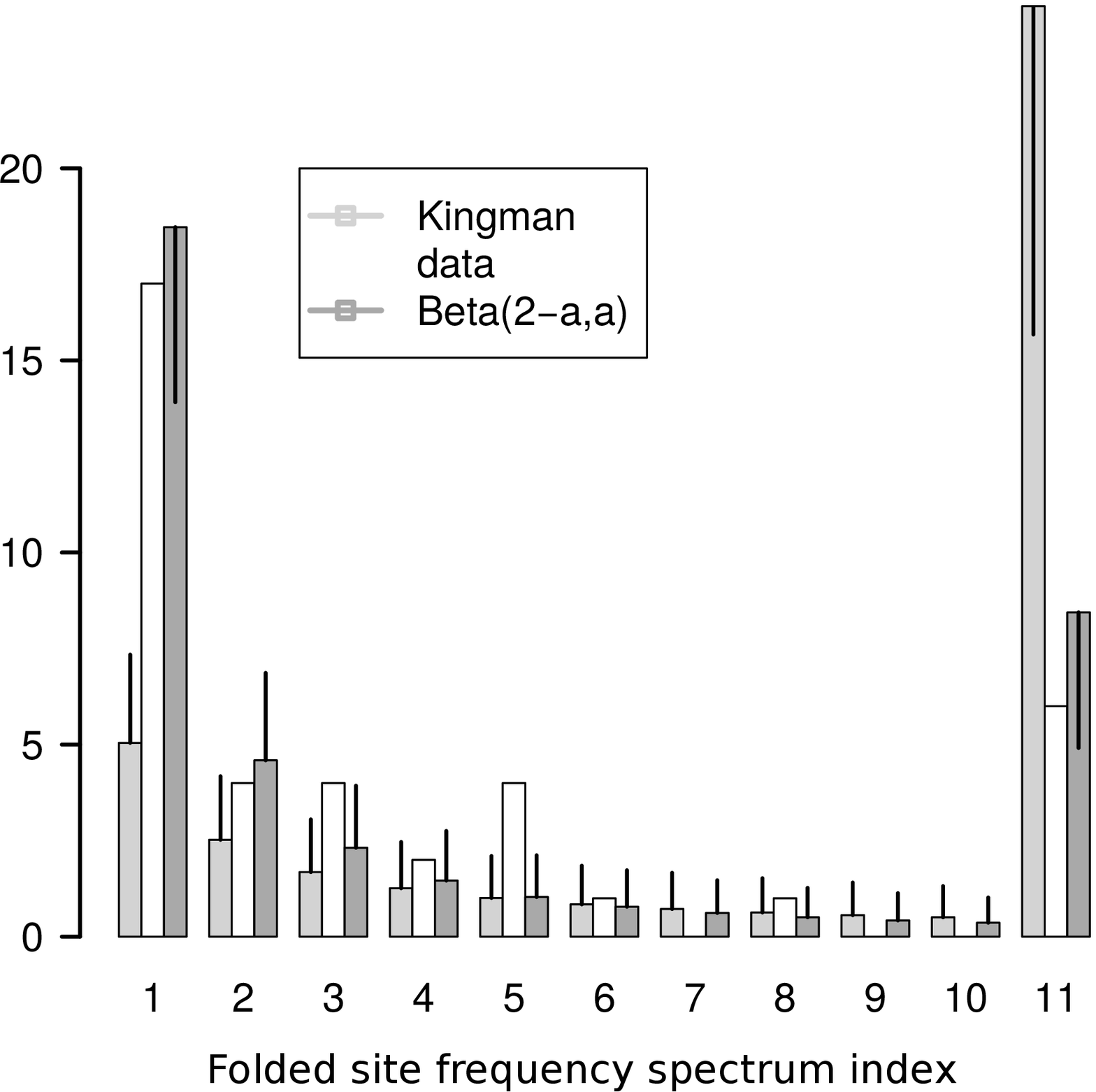} \hspace{0.5cm} %
\parbox[b]{5cm}{\small The folded freq.\ spectrum (white bars) of the data of \cite{MBJB-A04} along with predictions of the Kingman coalescent (light-grey), and the Beta$(2-\hat{\alpha},\hat{\alpha})$-coalescent (dark-grey), where $\hat{\alpha}=1.5$ is the best fit estimated from the data according to \cite{MBJB-BBE_SFS_13}. Vertical lines represent the standard deviation; obtained for the Beta$(2-\hat{\alpha},\hat{\alpha})$-coalescent from $10^5$ iterations. Class `11' represents the collated tail of the spectrum, from 11 to 1278/2.\\[1ex] Reproduced from \cite[Fig.~11]{MBJB-BBE_SFS_13}.
\\[0.1cm]}
\caption{} \label{MBJB-fig:codalldata}%
\vspace{-0.3cm}
\end{figure}

It is interesting to see that the expected site-frequency spectra differ significantly for the various coalescent models. In Figure~\ref{MBJB-fig:codalldata}, we compare the folded expected frequency spectra of a Kingman and a Beta-coalescent. We also include the frequency spectrum of mtDNA data for Atlantic cod from \cite{MBJB-A04} (1278 sequences). The fit of the Beta-coalescent to the real dataset is striking, see \cite{MBJB-BBE_SFS_13} for a discussion.

\begin{remark}
  1.\ For a $\Lambda$-coalescent $\Pi$ there are analogous recursions for variances 
  $\mathrm{Var}^{\Pi}\big[\xi^{(n)}_i\big]$ and covariances 
  $\mathrm{Cov}^{\Pi}\big[\xi^{(n)}_i, \xi^{(n)}_j\big]$, see
  \cite[Theorem~2]{MBJB-BBE_SFS_13}. \smallskip

  \noindent 2.\ For the Kingman case, we have
  ${p^{(n),{\delta_0}}[k,b]} = \frac{\binom{n - b - 1}{k - 2} }{\binom{n
      - 1 }{k - 1 }}$ and
  $\EE^{{\delta_0}, \theta} = \left[\xi_i^{(n)}\right] = \frac{\theta}{i}$,
  as computed by Fu \cite{MBJB-F95}. For general $\Lambda$-coalescents,
  no closed expressions for \eqref{MBJB-eq:pLambda_kb}, \eqref{MBJB-eq:Exii}
  are known. 
  However, the recursions can easily be solved numerically, even for $n$ in the hundreds.
  \smallskip

  \noindent 3.\ The computation of the expected SFS through \eqref{MBJB-eq:Exii} is natural and conceptually 
  appealing. We note however that there are now numerically more efficient alternatives, either via 
  a spectral decomposition of the jump rate matrix of $Y$ as in Spence~et~al \cite{MBJB-SKS16} or 
  via an interpretation as a multivariate phase-type distribution as in Hobolth~et~al's approach \cite{MBJB-HSJB18}.
  \smallskip

  \noindent 4.\ 
  For $\Lambda$-coalescents with `strong $\alpha$-regular variation' near $0$
  (i.e., $\Lambda(dx) = f(x)dx$ with $f(x) \sim A x^{1-\alpha}$ as $x\downarrow 0$ for some $A \in (0,\infty)$;
  this includes the $\mathrm{Beta}(2-\alpha,\alpha)$-coalescent 
  from Example~\ref{MBJB-ex:Lambdacoals}),
  \cite[Thm.~8]{MBJB-BBL14} shows $\xi^{(n)}_i \sim \frac{\theta}{2} n^{2-\alpha} C_{\alpha,i}$
  a.s.\ with an explicit constant $C_{\alpha,i}$. 
  However, the convergence in $n$ 
  can be quite slow,
  see \cite[Figure~8]{MBJB-BBE13} 
  and the discussion there.
  \smallskip

  \noindent 5.\ Using similar arguments, one can derive recursion formulas for the
  expectation and covariances of the site frequency spectrum under $\Xi$-coalescents.
  See 
  \cite{MBJB-BCEH16} and \cite{MBJB-SKS16}.
\end{remark}

\index{branch lengths of a coalescent tree}
We see from \eqref{MBJB-eq:trueLH} below and the following discussion that
the SFS is closely allied to the distribution of branch lengths in
coalescents.  Asymptotic results for such lengths are a focus of the
project by G.~Kersting and A.~Wakolbinger, described in this volume.  E.g.,
see \cite{MBJB-DK18a, MBJB-DK18b} for the asymptotic behaviour of
$B^{(n)}$ (the total branch length for sample size $n$) and of $B^{(n)}_1$ (the total branch length of the leaves) for very general coalescents
and \cite{MBJB-DKW14} for the fluctuations of $(B^{(n)}_1- \EE[B^{(n)}_1])/n^{1-\alpha+1/\alpha}$ for
$\mathrm{Beta}(2-\alpha,\alpha)$-coalescents with $1<\alpha<2$.
For the Bolthausen-Sznitman coalescent and some `relatives,'
corresponding to $\alpha=1$, \cite{MBJB-DK18a}
obtain the asymptotic behaviour as $n\to\infty$ of $B^{(n)}_i$ for any $i \in \NN$,
see 
the article by G\"{o}tz Kersting and Anton Wakolbinger in this volume.

The question of the theoretical identifiability of coalescents models
from the expected site frequency spectrum has been treated in
\cite{MBJB-SKS16}. For example for $\Lambda$-coalescents, the first $n-2$
moments of the measure $\Lambda$ can be determined from the expected
SFS with sample size $n$ and vice versa.

\subsection{Inference methods based on the site-frequency spectrum}
\label{MBJB-sect:InferenceSFS}

\subsubsection{Inference of mutation rates and real-time embeddings.}

When analysing data based on the SFS, one often needs to infer the underlying mutation rate first.
Hence we begin this subsection with 
a brief discussion of this estimation and 
its consequences for the real-time embedding (assuming a ``molecular clock'') of our coalescent models.
Estimating $\theta$ (or $\theta/2$) is 
often done via the (analogue of) the Watterson estimator. Here, as pointed out e.g.\ in \cite{MBJB-EBBF15}, it is important to understand that the choice of a multiple merger coalescent model $\Pi$
strongly affects this estimate. We illustrate this with an example. Assume w.l.o.g.\ for all multiple merger coalescents in
question that the underlying coalescent measure $\Lambda$ is always a
probability measure: This normalisation fixes the coalescent time 
unit as the expected time to the most recent common ancestor of two 
individuals sampled uniformly from the population. 

Given an observed number of segregating sites $S$
in a sample of size $n$, a common (and unbiased) estimate $\hat \theta^{\Pi}$ of the scaled mutation rate $\theta$ 
in the coalescent scenario $\Pi$ is the Watterson estimate
\begin{equation}
\label{MBJB-eq:thetaestimate}%
\hat \theta^{\Pi} := \frac{2S}{\EE^{\Pi}\big[B^{(n)}\big]}, 
\end{equation} 
where again $\EE^{\Pi}\big[B^{(n)}\big]$ is the expectation of the total tree length $B^{(n)}$
of an ($n$-) coalescent model $\Pi$. 
One can compute for example $\EE^{\Pi}\big[B^{(n)}\big]=\sum_{k=1}^n k g(n,k)$ with the Green function $g(n,k)$ 
from \eqref{MBJB-eq:YtGreenfct}.

Now with the estimate $\hat \theta^{\Pi}$, given knowledge of
the substitution rate $\hat \mu$ per year at the locus under consideration,
one can obtain an approximate real-time embedding of the coalescent history via 
\begin{equation}
\label{MBJB-eq:real_time_embedding}
{\rm coal.\ time\ unit} \times \frac{\hat \theta^{\Pi}}{2}\approx {\rm year} \times  \hat \mu. 
\end{equation}
cf.\ \cite[Section 4.2]{MBJB-SBB13}, which of course depends on the law $\PP^\Pi$ of the $\Pi$-coalescent via the expected value $\EE^{\Pi}\big[B^{(n)}\big]$. 
See also \cite{MBJB-WS09} for a study of the related concept of `effective population size.'

Given a Cannings population model of fixed size $N$ as discussed
in Section~\ref{MBJB-sect:popmodels}, let $c_N$
\index{pair coalescence probability}
be the probability that two gene copies, drawn uniformly at random and
without replacement from a population of size $N,$ derive from a
common parental gene copy in the previous generation. While for the
usual haploid Wright-Fisher model $c_N = 1/N$,
in the class \texttt{(B)} from Section~\ref{MBJB-sect:popmodels},
$c_N$ is proportional to
$1/N^{\alpha - 1}$, for $1 < \alpha \le 2$. By the limit theorem for
Cannings models of \cite{MBJB-MS2001}, one coalescent time unit
corresponds to approximately $1/c_N$ generations in the original model
with population size $N$. Thus the mutation rate $\tilde \mu$ at the
locus under consideration per individual per generation must be scaled
with $1/c_N$, and the relation between $\tilde \mu$, the coalescent
mutation rate $\theta^\Pi/2$ and $c_N$ is then given by the
(approximate) identity $c_N \approx 2 \tilde \mu/\theta^\Pi$.
In particular, if a Cannings model class (and thus $c_N$ as a
function of $N$) is given, 
the `effective population size' $N$ can then be estimated.
\index{effective population size}

\subsubsection{Approximate likelihood functions based on the SFS} 
\index{likelihood function}
\index{likelihood function!approximate}
Since mutations in our models occur as a Poisson process along the branches of a coalescent tree, 
for $\underline{k}=(k_1,k_2,\dots,k_{n-1})$ with $|\underline{k}|=\sum_{i=1}^{n-1} k_i =s$, the true likelihood function is 
\begin{align}
  \label{MBJB-eq:trueLH}
  L\left((\Pi, \theta), \underline{k}\right) 
  & = \PP^{\Pi,\theta} \big\{ \xi^{(n)}_i=k_i^{(n)}, i \in [n-1] \big\} 
  = \EE^{\Pi} \Bigg[ \prod_{i=1}^{n-1} e^{-\frac\theta2 B^{(n)}_i} \frac{(\theta B^{(n)}_i/2)^{k_i}}{k_i!} \Bigg] \notag \\
  & = \EE^{\Pi} \Bigg[ e^{-\theta B^{(n)}/2} \frac{(\theta B^{(n)}/2)^s}{s!} \cdot  
    \frac{s!}{k_1! \cdots k_{n-1}!}\prod_{i=1}^{n-1}\bigg(\frac{B_i^{(n)}}{B^{(n)}}\bigg)^{k_i}\Bigg]
\end{align}
where $B_i^{(n)}$ is the random length of branches subtending $i \in [n-1]$
leaves and $B^{(n)} = B_i^{(n)} + \cdots + B_{n-1}^{(n)}$ is the {\em total branch length} of the $n$-coalescent tree $\Pi$. 
\eqref{MBJB-eq:trueLH} is in general not expressible as a simple formula 
involving the coalescent parameters; it is in principle 
straightforwardly approximable via a `naive' Monte Carlo approach but this 
is computationally very expensive even for moderate sample sizes. 
We note that Sainudiin and V\'{e}ber \cite{MBJB-SV18} implement a clever approach to computing the expectation in \eqref{MBJB-eq:trueLH}
via importance sampling in the case of the Kingman coalescent (including variable population size and geographic 
structure); as far as we know, there is currently no study analogous to \cite{MBJB-SV18} that would include multiple 
merger coalescents.

Let us discuss an approximate likelihood function based on the
so-called `fixed-$s$-method'. The idea is to treat the observed number
of segregating sites as a {\em fixed parameter} $s \in \NN$, not as
(realisation of a) random variable $S$. 
This approximation appears quite common in the population genetics literature, 
see \cite{MBJB-EBBF15} and references there. 
Consider 
\begin{align}
\label{MBJB-eq:Likrepr}%
    \EE^\Pi \Bigg[ \frac{s!}{k_1^{(n)}! \cdots k_{n-1}^{(n)}!}\prod_{i=1}^{n-1}\bigg(\frac{B_i^{(n)}}{B^{(n)}}\bigg)^{k_i^{(n)}}\Bigg],
\end{align}
(i.e., we take only the last term inside the expectation in \eqref{MBJB-eq:trueLH}), this corresponds to 
uniformly and independently throwing $s$ mutations on the coalescent tree.
An approximation is 
\begin{equation} 
\label{MBJB-eq:Likrepr:proxy1} 
L(\Pi, \underline{k}^{(n)},s) \approx  
\frac{s!}{k_1^{(n)}! \cdots k_{n-1}^{(n)}!}\prod_{i=1}^{n-1}\big(\varphi^{\Pi,(n)}_i\big)^{k_i^{(n)}}
\end{equation}
where we replaced the random quantities $B_i^{(n)}/B^{(n)}$ in 
\eqref{MBJB-eq:Likrepr} by  the {\em expected normalised branch lengths} 
\begin{equation} 
\label{MBJB-eq:normSFS}
\varphi^{\Pi,(n)}_i = \EE^\Pi[B_i^{(n)}]/\EE^\Pi[B^{(n)}].
\end{equation}

Equation \eqref{MBJB-eq:Likrepr:proxy1} motivates the following family of 
`approximate' (in a twofold sense: regarding both fixing $s$ and exchanging expectation of a fraction with a fraction of expectations) likelihood functions
\begin{align}
\label{MBJB-eq:poissonfunction}
\widetilde L(\Pi, \underline{\xi}^{(n)};s) &=
\prod_{i=1}^{n-1}e^{- \frac{\hat \theta(\Pi, s)}2 \EE^\Pi[B^{(n)}] \varphi_i^{\Pi, (n)}}
\frac{\big(\frac{\hat \theta(\Pi, s)}2  \EE^\Pi[B^{(n)}] \varphi_i^{\Pi, (n)}\big)^{\xi^{(n)}_i}}{\xi^{(n)}_i!} \notag \\
&= \prod_{i=1}^{n-1} e^{-s \varphi_i^{\Pi, (n)}} \frac{(s \varphi_i^{\Pi, (n)})^{\xi^{(n)}_i}}{\xi^{(n)}_i!} 
\end{align} 
where $\hat\theta(\Pi, s)=2s/\EE^\Pi[B^{(n)}]$ is the Watterson estimator for the 
mutation rate under a $\Pi$-coalescent with $n$ leaves when $S=s$ 
segregating sites are observed, recall \eqref{MBJB-eq:thetaestimate}. 
In \eqref{MBJB-eq:poissonfunction}, we view $s$ as a parameter rather than 
as observed data, noting that $\widetilde L$ is well defined 
even if $|\underline{\xi} ^{(n)}| \neq s$. 

Note that for a principled approach to remove the dependence on the `nuisance parameter' $\theta$, 
one could follow \cite{MBJB-BB94}. However, this is computationally very costly in the 
context of MMC's and we do not pursue it here. For further discussion see \cite{MBJB-EBBF15}.
\smallskip

\eqref{MBJB-eq:poissonfunction} is a practical starting point for 
testing and parameter inference for multiple merger coalescent models, 
in particular this can be evaluated (and optimised) numerically very easily even for 
large sample sizes $n \gg 1$.

Let us also remark that \eqref{MBJB-eq:normSFS} can also be the starting point for 
inference based on minimum-distance statistics, see \cite{MBJB-BBE_SFS_13}.

\subsection{Can one distinguish population growth from multiple merger coalescents?}
\label{MBJB-subsect:Growth vs MMC}

\index{Beta-coalescent}
\index{Kingman coalescent!with exponential growth}
We now employ the approximate likelihood functions from the previous section to construct a likelihood-ratio test for model selection. While this method has also been employed to select between various $\Xi-$coalescent models (see \cite{MBJB-BBE13}), it can also be used to distinguish between different `evolutionary forces' leading to non-Kingman-like variability in the SFS.

As an example, we discuss a scenario where the underlying population in question has undergone an exponential population increase as in \cite{MBJB-EBBF15}. 
Consider a haploid Wright-Fisher model with population size $N$ at generation $r=0$ and size 
$N(r) = N(1 + \beta/N)^{-r}$ in generation $r$ before the present. This is in fact a special case of the set-up in \cite{MBJB-KK2003} 
and we obtain in the limit, by speeding up time with a factor $N$ as usual, a Kingman-coalescent with exponentially growing coalescence rates
$\nu(s)= e^{\beta s}$. Such a time-changed Kingman coalescent satisfies equation \eqref{MBJB-eq:masterSFS}.

A population which has undergone a recent rapid increase should
produce an excess of singletons in the SFS compared to model \texttt{(K)}, which is a pattern also
observed for Beta-coalescents. Similarly, Tajima's $D$ (a classical
test statistic in the Kingman context,
see \cite[Section~4.3]{MBJB-W08}) would tend to be significantly
negative under both model classes.

Our aim is to construct a statistical test to distinguish
between the model classes ${\tt(E)}$ and ${\tt(B)}$ (which intersect exactly
in ${\tt(K)}$). In order to distinguish ${\tt(E)}$ from ${\tt(B)}$, based on an observed
site-frequency spectrum $\underline{\xi}^{(n)}$ with sample size $n$ and
$S=|\underline{\xi}^{(n)}|$ segregating sites, a natural approach is to construct a  likelihood-ratio  test. 

Suppose our null-hypothesis $H_0$ is presence of recent exponential population
growth $({\scalebox{0.8}{\tt E}})$ with (unknown) parameter $\beta \in [0, \infty)$, and we
wish to test it against the alternative $H_1$ hypothesis of a multiple
merger coalescent, say, the Beta$(2-\alpha, \alpha)$-coalescent $(\scalebox{0.8}{\tt B})$ 
for  (unknown) $\alpha \in [1,2]$, where  $\beta = 0$ and $\alpha = 2$ 
correspond to the Kingman coalescent.
The coalescent mutation rate $\theta$ is not directly observable, 
but plays the role of a nuisance parameter.  By fixing $S=s$ and treating it as a parameter of our test, we 
may consider the pair of hypotheses
\begin{equation}
  H^s_0 \,:\, \Pi \in \Theta_s^{\scalebox{0.8}{\tt E}}:=\big\{ \text{Kingman coal., 
  growth parameter } \beta 
       : \beta \in [0,\infty)\big\} \vspace{-1ex}
\end{equation} 
and \vspace{-1ex}
\begin{equation}
  H^s_1 \,:\, \Pi \in \Theta_s^{\scalebox{0.8}{\tt B}}:=\big\{ \text{Beta$(2-\alpha, \alpha)$-coalescent} 
  : \alpha \in [1,2]\big\}. 
\end{equation} 
\index{likelihood-ratio text!approximate}
We can 
construct an  `approximate likelihood-ratio' test based on $L(\Pi, \underline{\xi}^{(n)},s)$
via 
\begin{equation}
\label{MBJB-eq:lira}
\varrho_{(\scalebox{0.8}{\tt E}, \scalebox{0.8}{\tt B};s)}(\underline{\xi}^{(n)}):= \frac{\sup\big\{L(\Pi, \underline{\xi}^{(n)},s),\,  \Pi \in \Theta_s^{\scalebox{0.8}{\tt E}}\big\}}{\sup\big\{L(\Pi, \underline{\xi}^{(n)},s),\, \Pi \in \Theta_s^{\scalebox{0.8}{\tt B}}\big\}}
\end{equation}
introduced in the previous section.
Given a significance level $a \in (0,1)$ (say, $a=0.05$), let $\varrho^*_{(\scalebox{0.8}{\tt E}, \scalebox{0.8}{\tt B};s)}(a)$ be the $a$-quantile of $\varrho_{(\scalebox{0.8}{\tt E}, \scalebox{0.8}{\tt B};s)}(\underline{\xi}^{(n)})$ under $\scalebox{0.8}{\tt E}$, chosen as the largest value so that 
\begin{equation}
\label{MBJB-eq:quantile}
\sup_{\Pi \in \Theta_s^{\scalebox{0.8}{\tt E}}} \PP^{\Pi,s} \big\{\varrho_{(\scalebox{0.8}{\tt E}, \scalebox{0.8}{\tt B};s)}(\underline{\xi}^{(n)}) \le \varrho^*_{(\scalebox{0.8}{\tt E}, \scalebox{0.8}{\tt B};s)}(a)\big\} \le a.
\end{equation}
The decision rule that constitutes the `fixed-$s$-likelihood-ratio test', given $s$ and sample size $n$, is
\vspace{-2ex}
$$
\hspace{0.5cm}\mbox{ reject $H_0^s$} \quad \iff \quad \varrho_{(\scalebox{0.8}{\tt E}, \scalebox{0.8}{\tt B};s)}(\underline{\xi}^{(n)}) \le \varrho^*_{(\scalebox{0.8}{\tt E}, \scalebox{0.8}{\tt B};s)}(a).
$$
The corresponding power function of the test, that is, the probability to reject a false null-hypothesis, 
is given by
\begin{equation}
\label{MBJB-eq:powerfunction}
G_{(\scalebox{0.8}{\tt E}, \scalebox{0.8}{\tt B};s)}(\Pi) = \PP^{\Pi} \{\varrho_{(\scalebox{0.8}{\tt E}, \scalebox{0.8}{\tt B};s)}(\underline{\xi}^{(n)}) \le
\varrho^*_{(\scalebox{0.8}{\tt E}, \scalebox{0.8}{\tt B};s)}(a)\}, \quad   \Pi \in \Theta_s^{\scalebox{0.8}{\tt B}}.
\end{equation} 
Alternatively, even though $\widetilde L(\cdot, \cdot \,;s)$ from \eqref{MBJB-eq:poissonfunction}
is not literally a likelihood 
function of any model from $H_0^s \cup H_1^s$, we can consider the statistic
$\widetilde\varrho_{(\scalebox{0.8}{\tt E}, \scalebox{0.8}{\tt B})}(\underline{\xi}^{(n)})$, where we replace
in \eqref{MBJB-eq:lira} $L(\Pi, \underline{\xi}^{(n)},s)$
by $\widetilde L(\Pi, \underline{\xi}^{(n)},|\underline{\xi}^{(n)}|)$.
For a given value of $s$, we can then (by simulations using the 
fixed-$s$-approach) determine approximate quantiles $\widetilde\varrho^*_{(\scalebox{0.8}{\tt E}, \scalebox{0.8}{\tt B};s)}(a)$
associated with a significance level $a$ as in \eqref{MBJB-eq:quantile}, 
and base our test 
on the criterion $\widetilde\varrho_{(\scalebox{0.8}{\tt E}, \scalebox{0.8}{\tt B})}(\underline{\xi}^{(n)}) 
\le \widetilde\varrho^*_{(\scalebox{0.8}{\tt E}, \scalebox{0.8}{\tt B};s)}(a)$.
Similarly, the (approximate) power function
\begin{equation}
  \label{MBJB-eq:tildeG_E_B}
  \widetilde G_{(\scalebox{0.8}{\tt E}, \scalebox{0.8}{\tt B};s)} = \PP^{\Pi} \{\widetilde \varrho_{(\scalebox{0.8}{\tt E}, \scalebox{0.8}{\tt B};s)}(\underline{\xi}^{(n)}) \le
  \widetilde \varrho^*_{(\scalebox{0.8}{\tt E}, \scalebox{0.8}{\tt B};s)}(a)\}
\end{equation} 
for $\Pi \in \Theta_s^{\scalebox{0.8}{\tt B}}$ can be estimated using simulations.
See the discussion in \cite{MBJB-EBBF15} and in particular Figure~2 there (a part of which we reproduce in 
Figure~\ref{MBJB-fig:newpower1} below).
For example, if the `truth' was a Beta$(2-\alpha, \alpha)$-coalescent with $\alpha=1.5$,
the power of a test of this form with significance level $5\%$ to reject $H_0^s$ 
(the null hypothesis of a Kingman model with exponential growth) 
based on a (single-locus) sample of size $n=500$ would
be about $75\%$. Note that the power is reasonably high for $\alpha \le 1.5$, say, but decays 
to the nominal level as $\alpha\to 2$. The boundary case $\alpha=2$ in the class of 
Beta$(2-\alpha, \alpha)$-coalescents \emph{is} the Kingman coalescent, after all.

\begin{figure}[h!]
  \ \hspace{-2cm}
  \begin{center}
    \includegraphics[width=7.0cm]{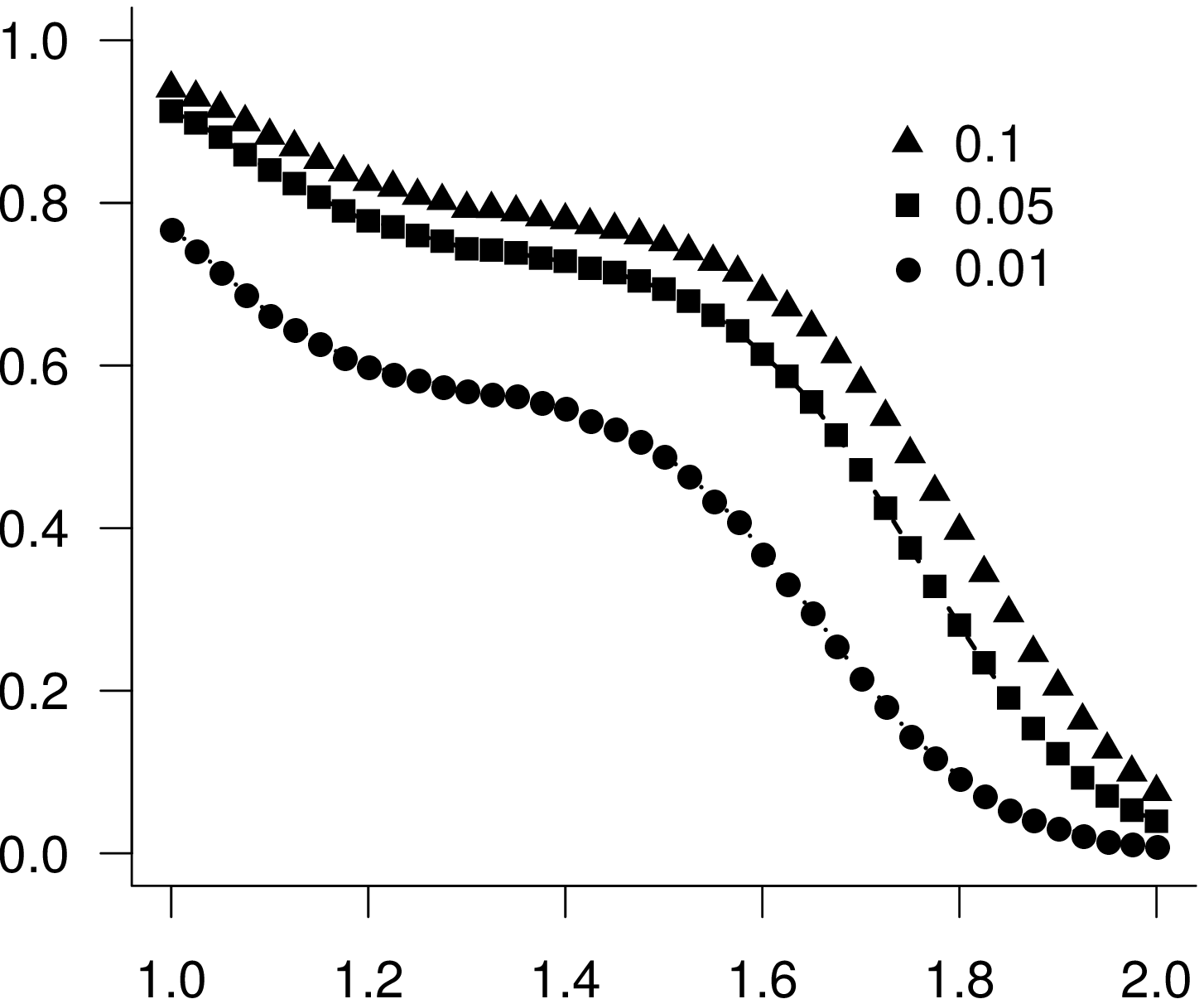}\\[0.5cm]
  \end{center}
  \caption{%
    Estimate of $\widetilde G_{(\protect\scalebox{0.8}{\tt E}, \protect\scalebox{0.8}{\tt B};s)}$ from \eqref{MBJB-eq:tildeG_E_B} 
    based on \eqref{MBJB-eq:poissonfunction} as a function of
    $\alpha$ with $n=500$ 
    and $s=50$. 
    The symbols denote the size of the test, cf. legend.  
    The hypotheses are discretised to $\Theta_s^{{\protect\scalebox{0.8}{\tt E}}}= \{\beta : \beta \in \{0,1,2, \ldots, 10, 20, \ldots, 1000\}\}$ and
    $\Theta_s^{\protect\scalebox{0.8}{\tt B}}= \{ \alpha : \alpha \in \{ 1, 1.025, \ldots, 2\}\}$. 
    Here, the Beta$(2-\alpha,\alpha)$-coalescent is the alternative. 
    Image reproduced from \cite[Figure~2]{MBJB-EBBF15}.
  }
  \label{MBJB-fig:newpower1}%
\end{figure}

\section{Multiple loci, diploidy and $\Xi$-coalescents}
\label{MBJB-sect:multilocus etc}
  
\subsection{A diploid  bi-parental multi-locus model}
\label{MBJB-subsect:diploidmodel}
\index{diploid}
\index{multi-locus model}

\begin{figure}
  \begin{center}
    \begin{tabular}{c}
      \fbox{\raisebox{2ex}{$\begin{array}{c} \mbox{\small time} \\ \uparrow \end{array}$} 
      \includegraphics[height=0.8cm, width=7.5cm]{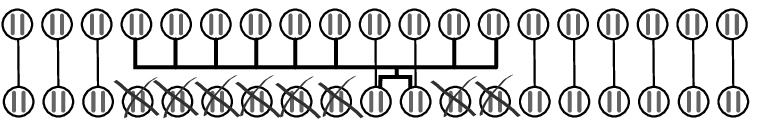}}
      \\[0.1cm]
      \fbox{\includegraphics[height=0.8cm, width=5.5cm]{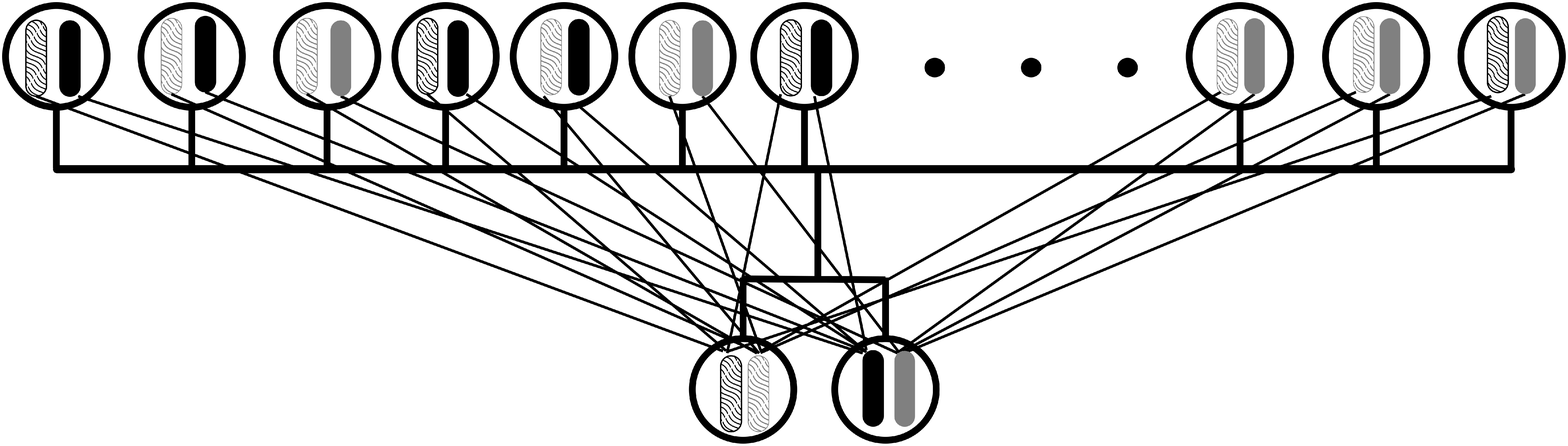}}
    \end{tabular}%
    \raisebox{-0.8cm}{\fbox{\includegraphics[height=1.9cm, width=1.2cm]{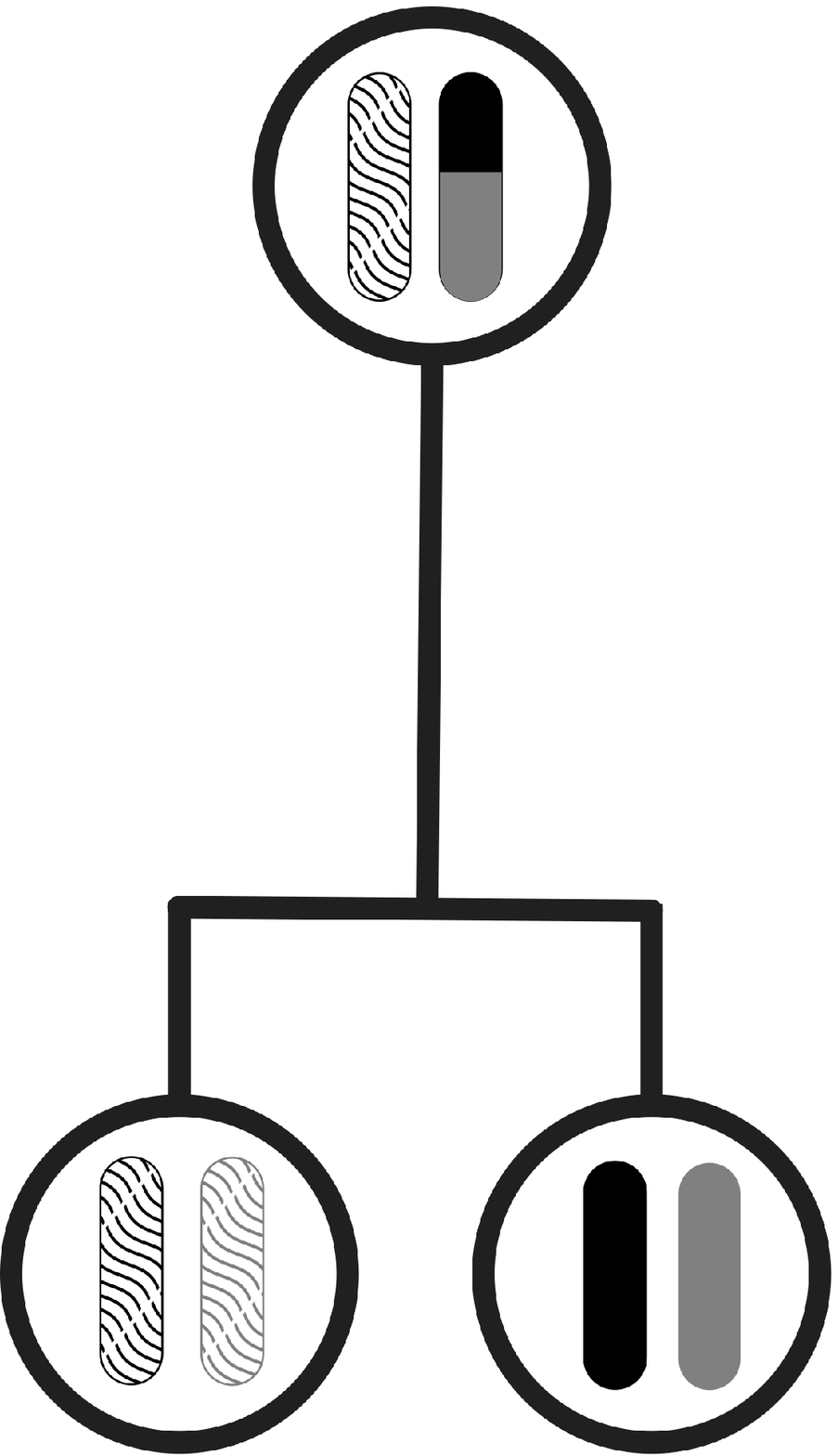}} 
      \fbox{\includegraphics[height=1.9cm, width=1.7cm]{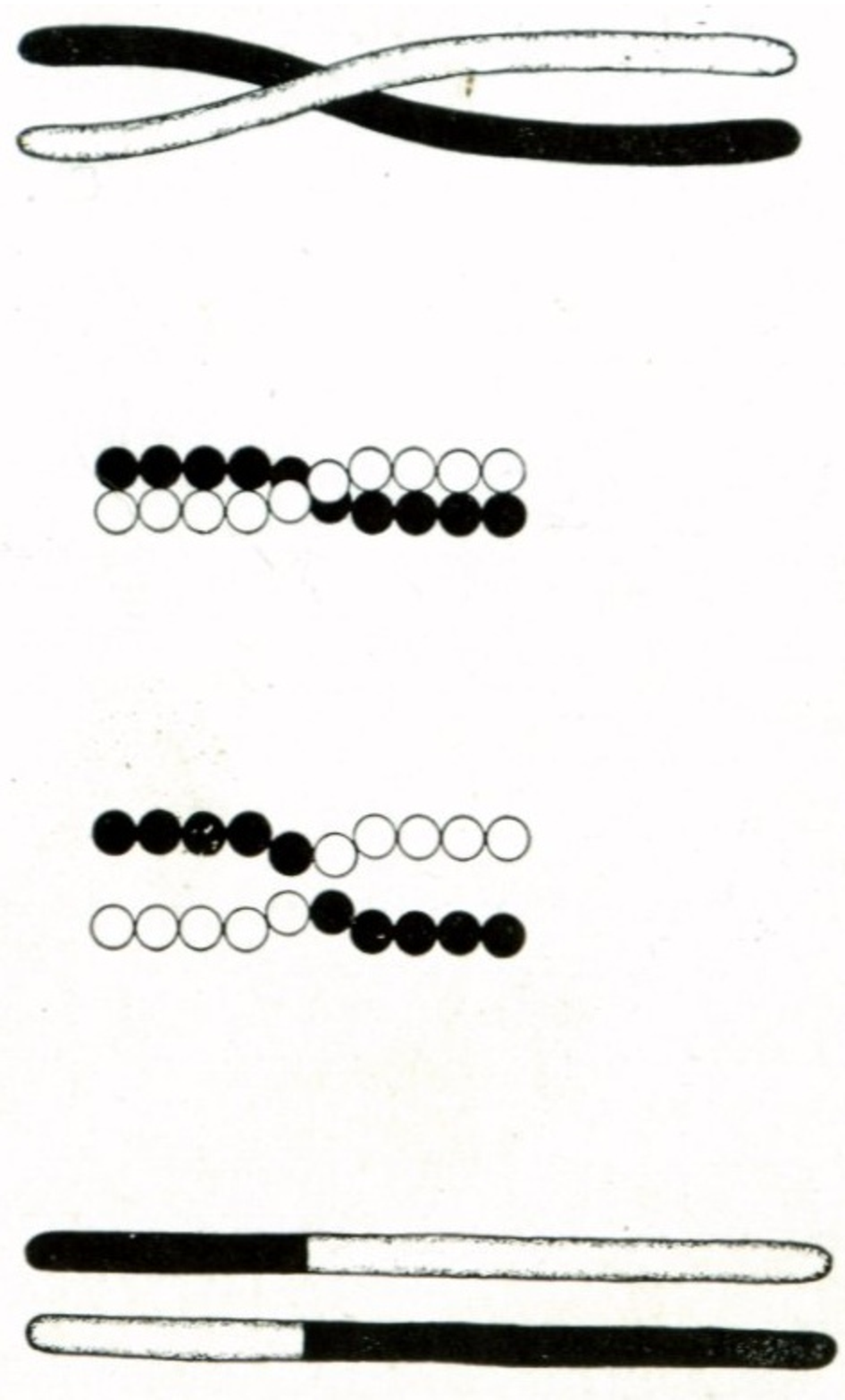}}}
  \end{center}
\caption{Schematic illustrations of the population model described in Section~\ref{MBJB-subsect:diploidmodel}. 
Top left: $\Psi^{(N)}$ children of a single pair are created. 
Bottom left: Transmission of genetic information to the $\Psi^{(N)}$ children (which can include recombination). 
Right: A possible recombination event in 
producing a child. 
Far right: Schematic illustration of crossing over (an important step in the biochemical mechanism of recombination), 
adapted from Thomas Hunt Morgan, \textit{A Critique of the Theory of Evolution}, Princeton University Press, 1916. 
}
\label{MBJB-fig:biparentmodel illustration}
\end{figure}

We model a population of $N$ diploid individuals. Each carries two chromosome copies, and each
chromosome consists of $L$ loci.  In a reproduction event, two
randomly chosen parents produce a random number $\Psi^{(N)}$ of
offspring, and these replace as many randomly chosen individuals;
$\Psi^{(N)}$ is drawn afresh for each event. 
Each child inherits one (possibly recombined) chromosome from each parent according to 
the Mendelian laws; we assume that during meiosis, a crossover recombination between 
locus $\ell$ and $\ell+1$ happens with probability $r_\ell^{(N)}$ for $\ell \in [L-1]$.
See Figure~\ref{MBJB-fig:biparentmodel illustration} for an illustration.

\begin{example}
For a concrete example, assume that $\PP(\Psi^{(N)} = \lceil \psi N \rceil) = c/N^2$ and 
$\PP(\Psi^{(N)} = 1) = 1-c/N^2$ with $\psi \in (0,1)$, $c>0$. This leads to model \texttt{(xEW)}.
\end{example}
\smallskip

Let $c_N := \EE\left[\Psi^{(N)}(\Psi^{(N)} + 3)/N(N - 1)\right]$
\index{pair coalescence probability}
(this $4 \times$ the pair coalescence probability for two randomly chosen chromosomes)
and assume that
\begin{align} 
  \label{MBJB-eq:cond time scale sep}
  \frac{c_N}{\EE\left[ \Psi^{(N)}/N\right]} 
  = \frac{\EE\left[\Psi^{(N)}(\Psi^{(N)} + 3)\right]}{(N-1) \EE\left[\Psi^{(N)}\right]}
    \mathop{\longrightarrow}_{N\to\infty} 0
\end{align}
(which implies that also $c_N\to 0$)
and that there exists a probability measure $\Lambda$ on $[0,1]$ such that
\begin{equation}%
  \label{MBJB-cond2:sagitovs}
  \frac{1}{c_N}\PP\left\{ \Psi^{(N)} > Nx \right\} 
  \mathop{\longrightarrow}_{N\to\infty} \int_{(x,1]} \frac{1}{y^2} \Lambda(dy) 
\end{equation}%
for all continuity points $x \in (0,1]$ of $\Lambda$. 
Furthermore 
\begin{equation} 
  \label{MBJB-eq:recombratescaling.general}
  r_\ell^{(N)} \sim \frac{c_N}{4\EE\left[\Psi^{(N)}/N\right]} r^{(\ell)}
  \quad \text{as}\; N\to\infty
\end{equation}
with fixed $r_\ell \in [0,\infty)$ for $\ell=1,\dots,L-1$.

\begin{remark}
\label{MBJB-rem:time scale Xi-ARG}
Note that $\EE\left[\Psi^{(N)} / N\right]$ is the probability that (after a
given reproduction event) a randomly chosen individual from the
current population is a child.
\eqref{MBJB-eq:cond time scale sep} then ensures that `separation of time
scales' occurs: The `short' time-scale $1/\EE\left[\Psi^{(N)} / N\right]$ on
which sampled chromosomes paired in the same individual disperse into 
two different individuals carrying only one sampled chromosome each is much smaller than the `long' time-scale $1/c_N$
over which we observe non-trivial ancestral coalescences. This lies 
`behind' Proposition~\ref{MBJB-thm:conv.generalXiARG} below.
\end{remark}

For the classification of general diploid models (in the single-locus context), we refer to \cite{MBJB-BLS18}, 
see also the article by Anja Sturm in this volume.

\subsection{The $\Xi$-ancestral recombination graph}
\label{MBJB-subsect:XiARG}

\index{ancestral recombination graph}
\index{$\Xi$-ancestral recombination graph}

Consider a sample of $n$ chromosomes (which could be taken from $n/2$ sampled individuals, say), each of 
which carries $L$ loci. 
We need some notation to describe the ancestral states: 
A possible 
configuration has the form $\zeta = \{ C_1, C_2, \dots, C_b \}$
with $b \in [n]$,
where $C_i = (\widetilde{C}_{i,1}, \widetilde{C}_{i,2},\dots, \widetilde{C}_{i,L})$ with 
$\widetilde{C}_{i,1}, \dots, \widetilde{C}_{i,L} \subset [n]$ and not all $= \emptyset$ such that
for $\ell=1,\dots,L$ we have $\bigcup_{i=1}^b \widetilde{C}_{i,\ell} = [n]$ 
and for $i \neq i'$, $\widetilde{C}_{i,\ell} \cap \widetilde{C}_{i',\ell} = \emptyset.$
$\widetilde{C}_{i,\ell}$ contains the indices of those samples for which the 
chromosome $C_i$ in the current configuration is ancestral at the $\ell$-th locus. 
Thus, for each locus $\ell$, $\{\widetilde{C}_{1,\ell},\dots,\widetilde{C}_{b,\ell}\}$ is a partition of 
$[n]$ (with a grain of salt: it may contain $\emptyset$'s).
We write $\cA$ for the set of all configurations of 
this form. 
We remark that in order to properly describe the dynamics of ancestral configurations for 
finite population size $N$, $\mathcal{A}$ 
is in fact not completely sufficient and has to 
be `enriched' by information about the grouping of ancestral chromosomes into diploid individuals.
However, because of the separation of time scales described in Remark~\ref{MBJB-rem:time scale Xi-ARG}, 
this becomes irrelevant for the limit process. We will not go into details here and 
refer to \cite{MBJB-BBE13}.

From $\zeta \in \cA$, 
possible transitions lead to \vspace{-1.0ex}
\begin{align*}
    \mathsf{pairmerge}_{i_1,i_2}(\zeta) & = \big\{
      C_1,\dots,C_{i_1-1},\widehat{C}_{i_1}, C_{i_1+1},\dots,
      C_{i_2-1},C_{i_2+1},\dots,C_b \big\} \\[-1.5ex]
    \intertext{with $\widehat{C}_{i_1} = \left(\widetilde{C}_{i_1,1} \cup \widetilde{C}_{i_2,1},\dots, 
  \widetilde{C}_{i_1,\ell} \cup \widetilde{C}_{i_2,\ell}\right)$, a merger of the pair $C_{i_1}$ and $C_{i_2}$,
  \vspace{-2.5ex}}
    \mathsf{groupmerge}_{J}(\zeta) & = \big\{ \overline{C}_1, \overline{C}_2, \overline{C}_3, \overline{C}_4, \,
                                     C_j, j \in [b] \setminus (J_1 \cup J_2 \cup J_3 \cup J_4 ) \big\} \\[-1.5ex]
  \intertext{with $J_1,\dots,J_4 \subset [b]$ pairwise disjoint and at least one $|J_i| \ge 3$ or 
  at least two of the $|J_i| \ge 2$. Here, $\overline{C}_m = \Big( \bigcup_{i \in J_m} \widetilde{C}_{i,1}, 
  \bigcup_{i \in J_m} \widetilde{C}_{i,2}, \dots, \bigcup_{i \in J_m} \widetilde{C}_{1,\ell}\Big)$ for $m=1,2,3,4$,
  a simultaneous multiple merger in (up to) four groups, and \vspace{-1.5ex}}
    \mathsf{recomb}_{i,\ell}(\zeta) & = \left\{ C_1,\dots,C_{i-1}, C_i', C_i'', C_{i+1},\dots,C_b \right\} 
\end{align*}
with
$C_i'=(\widetilde{C}_{i,1},\widetilde{C}_{i,2},\dots,\widetilde{C}_{i,\ell},\emptyset,\dots,\emptyset)$
and
$C_i'=(\emptyset,\dots,\emptyset,\widetilde{C}_{i,\ell+1},\widetilde{C}_{i,\ell+2},\dots,\widetilde{C}_{i,L},)$,
a recombination event splitting the $i$-th chromosome in the
configuration between locus $\ell$ and locus $\ell+1$.

Note that as mentioned above, both in the $\mathsf{pairmerge}$ and the
$\mathsf{groupmerge}$ operations, `empty' entries
$(\emptyset, \emptyset, \dots, \emptyset)$ may arise, which then need
to be removed; see \cite{MBJB-BBE13} for details.
\smallskip

\begin{figure}

  \includegraphics[width=6.5cm]{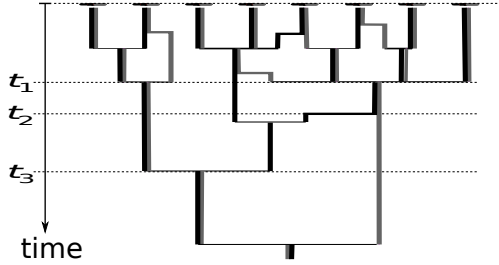} \ \vspace{0.1cm}
\caption{An illustration of the $\Xi$-ancestral recombination graph for two loci, with 
some transitions highlighted. At time $t_1$, a $\mathsf{groupmerge}$-event occurs. At time $t_2$, there is 
a $\mathsf{recomb}$-event and at time $t_3$, a $\mathsf{pairmerge}$-event.}
\label{MBJB-fig:XiARG0}
\end{figure}

The limiting genealogical process will then be a continuous-time
Markov chain $\{ \xi(t)\}_{t \ge 0}$ on $\cA$ with generator matrix $q$ whose off-diagonal
elements are given by 
  \begin{align}
    \label{MBJB-eq:DefG.general}
    q_{\xi,\xi'} = 
    \begin{cases} 
      C_{\beta;2}
      & \text{if} \; 
      \xi'= \mathsf{pairmerge}_{j_1,j_2}(\xi) \\
      r^{(\ell)} & \text{if} \; \xi'= \mathsf{recomb}_{j,\ell}(\xi) \\
      C_{\beta;|J| } 
      & \text{if} \; \xi'= \mathsf{groupmerge}_{J_1,J_2,J_3,J_4}(\xi) \\
      0 & \text{for all other $\xi'\neq\xi$}
    \end{cases}
  \end{align}
  where $C_{\beta;|J|} :=
  C_{\beta;|J_1|,|J_2|,|J_3|,|J_4|;\beta-(|J_1|+|J_2|+|J_3|+|J_4|)}$ and 
  \begingroup
  \allowdisplaybreaks
  \begin{align} 
    \label{MBJB-eq:Xirate.general}
    C_{b;k;s} 
   & = \Lambda(\{0\}) \delta_{\{r=1, k_1=2\}} 
     + 4 \sum_{l=0}^{s \wedge (4-r)} {s \choose l}  \tfrac{(4)_{r+l}}{4^{|k| + l}} 
      \cdot  \int_{(0,1]} x^{|k|+l} (1 - x)^{s - l} \frac{1}{x^2} \Lambda(dx) 
  \end{align}
  \endgroup
  with $k = (k_1, \ldots, k_r)$, $|k| = k_1 + \cdots + k_r$.
  The path of $\{\xi_t\}$ 
  can be visualised as a random network, see Figure~\ref{MBJB-fig:XiARG0} for an illustration.

\begin{proposition}[{\cite[Theorem~1.3]{MBJB-BBE13}}] 
  \label{MBJB-thm:conv.generalXiARG}
  Let $\{\xi^{n,N}(m), m \ge 0\}$ be the ancestral process of a sample of $n$ chromosomes
  in a population of size $N$ with offspring laws $\mathcal{L}(\Psi^{(N)})$ 
  satisfying 
  \eqref{MBJB-eq:cond time scale sep}
  and \eqref{MBJB-cond2:sagitovs}, and 
  assume the scaling relation \eqref{MBJB-eq:recombratescaling.general}.
  \begin{equation}
    \label{MBJB-thm:conv.generalXiARG_claim}
    \{\xi^{n,N}(\lfloor 4t/c_N \rfloor)\} \longrightarrow \{\xi(t)\} \quad \mbox{ as } \quad N \to \infty ,
  \end{equation}
  where the process
  $\{\xi(t)\}$ is the Markov chain with generator matrix 
  \eqref{MBJB-eq:DefG.general}.
\end{proposition}
We refer to \cite{MBJB-BBE13} for details, in particular the precise mode of convergence in \eqref{MBJB-thm:conv.generalXiARG_claim}
depending on whether or not the grouping of ancestral chromosomes into possibly `doubly marked individuals' 
is taken into account.

\subsection{Towards a full {SMMC} multilocus inference machinery}
\label{MBJB-section full multilocus}

One can incorporate the (biologically important) effects of
recombination, spatial subdivision, variable population size (e.g.\
growing populations), and/or (directional) selection into stochastic
models for populations with highly skewed offspring distributions and
derive corresponding (limiting) models for the joint genealogy of an
$n$-sample observed at $L$ (possibly recombining) loci. The `full
complexity' model is then a `structured $\Xi$-ancestral selection
recombination graph.'  While in principle highly relevant in view of
today's large scale datasets, an explicit description of the resulting
full sampling distributions 
seems out of reach at present. One can however
make progress on statistical questions by employing low-dimensional
summary statistics. One approach, inspired by the results from
Section~\ref{MBJB-sect:InferenceSFS} is to use suitable lumpings of
the normalised site frequency spectra and average these over the
observed loci: Let
\begin{align} 
  \zeta_1(\ell) := \frac{\xi_1(\ell)}{|\xi(\ell)|}, \quad 
  \overline\zeta_k(\ell) := \sum_{j=k}^{n-1} \frac{\xi_j(\ell)}{|\xi(\ell)|}
\end{align}
be the proportion of singletons and the proportion of mutations visible in more than $k \ge 2$ 
copies at the $\ell$-th locus, respectively. 
\begin{align}
  \label{MBJB:zeta2}
  \left( \zeta_1, \overline\zeta_k \right) 
  := \frac{1}{L} \sum_{\ell=1}^L \left( \zeta_1(\ell), \overline\zeta_k(\ell)\right)
\end{align}
is a two-dimensional summary of the data whose distribution under a given coalescent model 
$\Pi$ with mutation parameter $\theta>0$ 
\begin{align}
  \label{MBJB:2dsumlik}
  L\left(\Pi,\theta, (z_1, \overline{z}_k)\right) 
  := \PP^{\Pi,\theta}\Big( \left( \zeta_1, \overline\zeta_k \right) = (z_1, \overline{z}_k) \Big)
\end{align}
is generally not known explicitly, but $\left( \zeta_1, \overline\zeta_k \right)$ can be simulated 
readily under $(\Pi, \theta)$. Then the function $(z_1, \overline{z}_k) \mapsto 
L\left(\Pi,\theta, (z_1, \overline{z}_k)\right)$
from \eqref{MBJB:2dsumlik} can be approximated by a kernel estimator based on 
$M$ independent replicates: 
\begin{align}
  \label{MBJB:2dsumlikappr}
  \widehat{L}\left(\Pi,\theta, (z_1, \overline{z}_k)\right) 
  := \frac{1}{M h} \sum_{m=1}^M K\left( \frac1h \left( (\zeta_1, \overline\zeta_k) 
  - \big( \zeta_1, \overline\zeta_k \big)(m) \right) \right)
\end{align}
where $\big( \zeta_1, \overline\zeta_k \big)(m)$ is the value of
\eqref{MBJB:zeta2} computed from the $m$-th simulation and $K$ the
kernel function (e.g.\ a Gaussian) with bandwidth $h>0$.  Given
\eqref{MBJB:2dsumlikappr}, testing and model selection analogous to
Section~\ref{MBJB-subsect:Growth vs MMC} can now be based on the
approximate likelihood ratio statistic
\begin{align} 
  \label{MBJB:eqLhatratiostatistic}
  \frac{\sup_{(\Pi,\theta) \in \Theta_0} \widehat{L}\left(\Pi,\theta, (z_1, \overline{z}_k)\right)}{%
  \sup_{(\Pi,\theta) \in \Theta_1} \widehat{L}\left(\Pi,\theta, (z_1, \overline{z}_k)\right)}
\end{align}
where of course the critical value for a test of given size has to be 
determined by simulations. In practice, one can alleviate the two-dimensional optimisation problem 
in \eqref{MBJB:eqLhatratiostatistic} by plugging in the Watterson estimator $\theta = \hat{\theta}^\Pi$ 
from \eqref{MBJB-eq:thetaestimate} given coalescent model $\Pi$.

This approach is pursued in \cite{MBJB-K18}, 
with promising initial results, see the discussion there and also 
Figure~\ref{MBJB-fig:2dmultilocussummary} below. It can also be extended to include the effects of
selection, variable population sizes and spatial structure, see \cite{MBJB-KWB19}
for steps in this direction. Note that this is akin to approximate Bayesian computations (ABC), whose r\^{o}le in 
analyses of datasets in multiple merger contexts is described in the article by Fabian Freund 
in this volume. 

Intuitively, although
even unlinked loci are not independent under the skewed offspring
distribution models from Section~\ref{MBJB-subsect:XiARG} (as observed
in \cite{MBJB-BBE13}), averaging over many loci does reduce sampling
variability and is justified because the multiple merger mechanism
affects all loci in the same way.  This is in fact a distinguishing
feature that explains why multi-locus data is useful to distinguish
skewed offspring distributions from selective sweeps: The latter would
only affect one locus at a time.
\begin{figure}[h!]
  \begin{center}
    \includegraphics[scale=0.7]{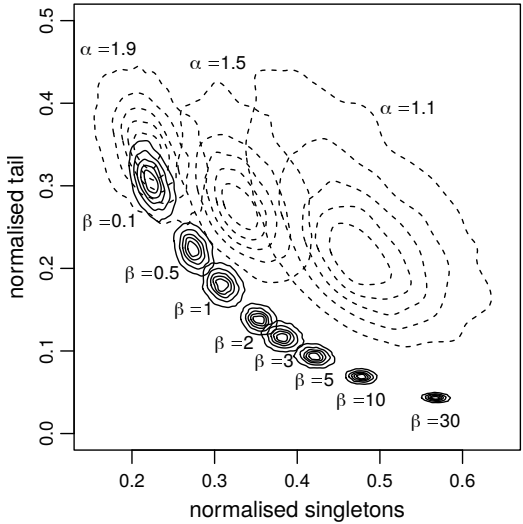}
  \end{center}
  \caption{The empirical distribution of $(\zeta_1, \overline\zeta_k)$ from \eqref{MBJB:zeta2}
  is quite different under a Kingman coalescent with exponential growth (solid contours) compared to a
  4-fold Beta coalescent \texttt{(xB)} (dashed contours).
  Here, the sample size is $n=100$, each sample considered at $L=23$ loci, with cutoff parameter $k=15$.
  Parameter values ($\alpha$ for the 4-fold Beta($2-\alpha,\alpha$) coalescent, $\beta$ for the
  exponential growth rate) are as shown. 
  The contour lines are based on 5000 simulated replicates
  for each parameter choice: For this, mutation rates $\theta$
  were chosen so that the expected number of segregating sites
  per locus equalled $s_\text{expect}^{(n),\Pi} = 10,20,30,40,50$ (cf. Equ.~\eqref{MBJB-eq:thetaestimate}), with 1000
  replicates per value of $\theta$. The pictures for a fixed value of $s_\text{expect}^{(n),\Pi}$
  are almost indistinguishable from the one shown.
  The contours were computed using \texttt{R} \cite{MBJB-R core team} and
  the function \texttt{kde} from the contributed R-package \texttt{ks} \cite{MBJB-Duong ks package},
  with default values for the bandwiths. 
  They correspond to regions containing respectively
  20\%, 40\%, 60\%, 80\% and 95\% of the simulated points.
  \label{MBJB-fig:2dmultilocussummary}}
\end{figure}

The software used for this study is available under \url{https://github.com/JereKoskela/Beta-Xi-Sim}. 
Furthermore, software for simulation and analysis of datasets in (S)MMC contexts can be found 
on Bjarki Eldon's homepage \url{http://page.math.tu-berlin.de/~eldon/programs.html}.

\section{Discussion - Are they really out there?}

In the previous sections, we outlined population models and
evolutionary scenarios which invite genealogical modelling via (S)MMC
processes. Further, we presented some paradigmatic statistical tools
for inference and model selection for (S)MMC processes, and our hope
is that this could pave at least some of the way towards an answer to
initial question \cite{MBJB-EW06} whether (S)MMC coalescents are really more
adequate null-models for real populations exhibiting highly skewed
offspring distributions (or other forces leading to an `effective
skew', such as selective sweeps, severe bottlenecks etc.).

One of our main take-home messages is that the statistical power of
such inference methods is usually much higher in (diploid) multi-locus
setups rather than in (haploid) single locus scenarios. However, it is
the latter scenario in which MMC based inference methods have so far
been applied in practice. For example, the results in
\cite{MBJB-SBB13} indicate that data generated under a
Beta-coalescent can provide a better fit to observed genetic
variability in Atlantic cod mitochondrial (thus haploid) DNA sequence
data. In the cited article, it is also discussed in how far 
different underlying coalescent models lead to different estimates for
the real-time most recent common ancestor of the sample.  
To some degree, it appears also possible to
distinguish different evolutionary scenarios such as a recent increase
in population size, leading to a time-changed Kingman coalescent, from
other coalescent scenarios, as reviewed in 
in Sections~\ref{MBJB-subsect:Growth vs MMC} and \ref{MBJB-section full multilocus}.

A very recent further study involving virus data (influenza) is
\cite{MBJB-SHJ19}, which employs purely-atomic MMCs (of class {\tt
  (EW)}), again in a haploid setup. The authors here come to the
conclusion that the {\tt (EW)} coalescent can provide a ``much more
accurate neutral null model'' in certain types of organisms including
viruses and bacteria. However, the study seems to be restricted to a relatively small class of MMCs.

We expect that a real test for the above methods will be in the framework of diploid multi-locus setups.
A very interesting step in this direction is the recent work of Rice, Novembre and Desai \cite{MBJB-RND18} 
who propose a statistic 
based on the joint site frequency spectrum at two loci. This approach does not explicitly 
model multi-locus dynamics including recombination, but it can (quite straightforwardly) be scaled up to 
analyse genome-wide genetic variability and, as shown in \cite{MBJB-RND18}, does shed a very interesting 
light on a Zambian population of fruit flies (\emph{Drosophila melanogaster}). 
Furthermore, in this context, it is rather satisfying to see that the funding of the Icelandic Grant of Excellence ``Population genomics of highly fecund codfish'' has recently been awarded jointly to \'{A}rnason, Halld\'{o}rsd\'{o}ttir, Etheridge, and  Stephan. Our hope is that this project will provide and analyse the necessary data on which the full multi-locus machinery can be tested. We will be curious to observe the outcomes.

\medskip

\noindent {\bf Acknowledgements.} The authors would like to thank
Iulia Dahmer, Frederik Klement and Timo Schl\"{u}ter for carefully reading the
manuscript and for their helpful comments. We also thank Iulia Dahmer for her help in
producing Figure~\ref{MBJB-fig:2dmultilocussummary} and two anonymous referees for their
insightful comments which helped to improve the presentation of this article.

\printindex

\end{document}